\documentclass{article}

\usepackage{subfiles}
\usepackage{amsmath,amssymb,bm,bbm,amsthm}
\usepackage{appendix}
\usepackage{lineno}
\usepackage{natbib}
\usepackage{tabularx}
\usepackage{graphicx, float}
\usepackage{geometry}
\usepackage{subcaption}
\usepackage{color}
\usepackage{makecell, booktabs}
\usepackage{caption}
\usepackage{url}

\newtheorem{proposition}{Proposition}
\newtheorem{definition}{Definition}

\let\originalleft\left
\let\originalright\right
\renewcommand{\left}{\mathopen{}\mathclose\bgroup\originalleft}
\renewcommand{\right}{\aftergroup\egroup\originalright}
\newcommand{\paren}[1]{\left(#1\right)}
\newcommand{\bracket}[1]{\left[#1\right]}
\newcommand{\braces}[1]{\left\{#1\right\}}

\newcommand{\Hstat}{G}
\newcommand{\pstat}{q}
\newcommand{\Dstat}{\mathcal{D}}

\newcommand{\st}{\text{s.t. }}

\newcommand{\celsius}{\text{C}}
\newcommand{\fahrenheit}{\text{F}}

\newcommand{\gens}{\mathcal{G}} 
\newcommand{\genscont}{\mathcal{G}^2} 
\newcommand{\gensbin}{\mathcal{G}^1} 
\newcommand{\N}{\mathcal{N}}
\newcommand{\None}{\mathcal{N}_1}
\newcommand{\Ntwo}{\mathcal{N}_2}

\newcommand{\sone}{x} 
\newcommand{\stwo}{y} 
\newcommand{\shed}{s} 

\newcommand{\timeperiods}{\mathcal{H}}
\newcommand{\Time}{\mathbf{H}}
\renewcommand{\time}{H}
\newcommand{\Temp}{\mathbf{T}}
\newcommand{\temp}{T}
\newcommand{\Outages}{\mathbf{Z}}
\newcommand{\outages}{Z}
\newcommand{\Demand}{\mathbf{D}}
\newcommand{\demand}{D}
\newcommand{\Rand}{{\bm \xi}}
\newcommand{\rand}{\xi}
\newcommand{\RandE}{\Rand_E}
\newcommand{\randE}{\rand_E}
\newcommand{\dist}{\mathbb{P}}
\newcommand{\distE}{\mathbb{P}_E}

\newcommand{\costcap}{c}
\newcommand{\costgen}{u}
\newcommand{\costgena}{u_N}

\newcommand{\costshed}{\lambda}
\renewcommand{\cap}{p}
\newcommand{\E}{\mathbb{E}} 
\renewcommand{\P}{\mathbb{P}} 
\newcommand{\cvar}[1][]{\text{CVaR}^#1} 
\newcommand{\chance}{\sigma}
\newcommand{\chanceparam}{v}
\newcommand{\cvarthres}{\alpha}
\newcommand{\cvarvar}{z}
\newcommand{\cvarvark}{v}
\newcommand{\lolp}{\text{LOLP}}
\newcommand{\lolpvark}{w}
\newcommand{\bigM}{M}
\newcommand{\obja}{Q_N}
\newcommand{\objb}{Q_E}
\newcommand{\recfeasa}{Y_N}
\newcommand{\recfeasb}{Y_E}
\newcommand{\recvara}{s_N}
\newcommand{\recvarb}{s_E}
\newcommand{\bound}{U}
\newcommand{\shedfn}{\ell_N} 
\newcommand{\indicator}{\mathbbm{1}}

\newcommand{\vone}{x}
\newcommand{\vtwo}{\bar{x}}

\newcommand{\xhat}{\hat{x}}
\newcommand{\ghat}{\hat{g}}

\begin{document}

\title{A Framework for Balancing Power Grid Efficiency and Risk with Bi-objective Stochastic Integer Optimization}
\author{Ramsey Rossmann \and Mihai Anitescu \and Julie Bessac \and Michael Ferris \and Mitchell Krock \and James Luedtke \and Line Roald}

\date{28 May 2024}

\maketitle 

\begin{abstract} 
    Power grid expansion planning requires making large investment decisions in the present that will impact the future cost and reliability of a system exposed to wide-ranging uncertainties.
    Extreme temperatures can pose significant challenges to providing power by increasing demand and decreasing supply and have contributed to recent major power outages.
    We propose to address a modeling challenge of such high-impact, low-frequency events with a bi-objective stochastic integer optimization model that finds solutions with different trade-offs between efficiency in normal conditions and risk to extreme events.
    We propose a conditional sampling approach paired with a risk measure to address the inherent challenge in approximating the risk of low-frequency events within a sampling based approach.
    We present a model for spatially correlated, county-specific temperatures and a method to generate both unconditional and conditionally extreme temperature samples from this model efficiently.
    These models are investigated within an extensive case study with realistic data that demonstrates the effectiveness of the bi-objective approach and the conditional sampling technique.
    We find that spatial correlations in the temperature samples are essential to finding good solutions and that modeling generator temperature dependence is an important consideration for finding efficient, low-risk solutions.
    \end{abstract}

\ifx\dissertationcompile\undefined
    \newcommand{\subdir}{}
\else
    \newcommand{\subdir}{ch1/}
\fi

\section{Introduction}
The problem of capacity expansion planning in the power grid involves making large capital decisions, like what generators to build and where to build them, in the face of an uncertain future. 
The decisions made now will impact the future cost and reliability of a system that will experience a wide range of uncertain loads and availability of capacity.
The system can also be significantly challenged by extreme events like storms and extreme temperatures that can have severe impacts on both supply and demand.
Many uncertainties in the power system are highly correlated due to their dependence on weather, especially temperature.
Extreme temperatures tend both to increase power demand and to decrease power supply \citep{murphy2019} and have contributed to major power outages and near-shortage events in recent years.
We present models and methods for supporting strategic decisions such as generator expansion in systems impacted by  high-impact, low-frequency events. 

The potential for severe consequences from such events, such as the widespread power outages caused by winter storms in Texas in February 2021 \citep{busby2021}, indicates a need for methods that allow decision-makers to balance the risk from these events with the everyday performance (efficiency) of the system. 
Given that protecting against all possible events would be prohibitively expensive, power system planners must find a trade-off between efficiency and risk. 
In order to model this trade-off, we propose a bi-objective stochastic integer optimization approach, where one objective is to minimize expected cost (including capital and operating costs), and the other is to minimize a measure of the risk of impacts from extreme events. 
A key feature of our approach is that the model for estimating impacts from extreme events can be different from the model for estimating costs from normal operation.
The goal of solving the bi-objective problem is to determine the set of Pareto optimal solutions, those solutions for which one objective cannot be improved upon without making the other objective worse. Hence, this yields solutions with different trade-offs between cost and risk. 
This allows system planners to make informed decisions that balance these competing objectives.

Our proposed models involve optimizing expected value (for system efficiency in normal operations) or a risk measure (for the risk from extreme events), which is generally intractable to do directly due to the challenge in evaluating the expected value or the risk measure over a general distribution. Sample average approximation (SAA) addresses this challenge by replacing the expected value (or the risk measure) with an empirical estimate obtained from a random sample of the uncertain quantities. However, in our context of systems impacted by rare, high-impact events, the sample size required to obtain good solutions may be intractably large since low probability events are then important for decision making.
As a result, modeling the risk to rare, high-impact events for the power grid can encounter computational challenges.
To address this, we propose a conditional sampling approach designed to obtain good quality solutions with modest sample sizes.
Given a distribution of uncertain parameters impacting the system, the idea is to use problem insight to determine characteristics of the parameters that are likely to result in bad outcomes for the system, and then sample from the distribution conditional on those characteristics. 
These conditional samples can then be used in SAA to approximate the risk objective using a relatively smaller sample size. We interpret this as a method for approximating the set of Pareto optimal solutions in a bi-objective optimization problem by replacing one objective with an approximation that is easier to optimize. In this context, we investigate the potential of this approximation by providing a sufficient condition for such an approach to yield an approximation of the set of Pareto optimal solutions.

Given the influence of temperature on the power grid, we pay special attention to how the distribution of temperatures is modeled.
As shown in \citet{ieee2019}, electric demand depends on temperature, due largely to air conditioner usage in hot weather and electric heat usage in cold weather. 
\citet{murphy2019} found generator availability to be temperature dependent, with colder and hotter temperatures leading to higher rates of generator outages. 
With the increasing demand in these situations as well, there is the potential for extreme temperatures to pose significant challenges for the power grid, as was the case in the 2021 Texas Power Crisis. 
As extreme temperature events (among other severe weather events) can be very widespread geographically, the ability of neighboring regions to provide power to an affected region is an important consideration. 
This makes spatial correlation in temperature distributions important for accurately assessing the system-wide impacts of extreme events.
We advocate using spatially correlated temperature models, and we propose a modeling approach and conditional sampling method for the case of extreme temperatures across a large region.

We conduct an extensive case study with realistic data to investigate the potential of the proposed modeling framework.
The case study considers generator investment decisions for the Midwestern United States, with temperature dependence for demand as well as generator capacities modeled on a county-by-county basis. 
Historical temperatures and actual generator specifications are used in modeling the costs and temperature effects using data from NREL and the EPA and findings from \citet{murphy2019}. 
Solar and wind capacity factors are modeled at the county level with an hourly resolution using data from NREL.
Results from the case study demonstrate that the conditional sampling used in the bi-objective approach is effective at finding Pareto optimal solutions at varying levels of risk and cost. 
We also find that modeling spatial correlations in the temperature samples is essential to finding good solutions and that the generator temperature dependence modeling is an important consideration for finding efficient, low-risk solutions.

\subsection{Relationship to existing work}
Techniques for optimization under uncertainty can be powerful tools for decision-making in power system planning (see \citet{roald2023} for a survey of power systems applications). 
Generation and transmission expansion have been modeled with stochastic programming considering uncertainty in demand to inform capacity expansion (\citet{delgado2013} and \citet{alvarez2006}) and unit commitment \citep{zheng2014}. 
Supply uncertainties, especially power generation from wind and solar resources, have also been modeled with stochastic programming to inform investment and dispatch decisions \citep{papavasiliou2011,banzo2011,cory2020stochastic}. 
Our work builds on existing capacity expansion models by explicitly allowing two modes of operation (``normal'' and ``extreme'' in our case) with different objectives. 
While we use simple operational subproblems in our case study, the flexibility of our bi-objective framework allows decision makers to tailor the model to different operational considerations and uncertainties.

Solving stochastic programs with sampling techniques is generally effective, and various approaches exist to address sampling issues of severe tail events \citep{ruszczynski2006}. 
Importance sampling \citep{tokdar2010} can reduce the required sample size by oversampling tail scenarios (then correcting for this bias a posteriori); however, doing so in stochastic programming has a circularity challenge since the magnitude of the impact from random events depends on the decisions made \citep{ruszczynski2006} and hence the best set of ``tail scenarios'' to oversample from are not known when performing the sampling. 
While prior work has attempted to address this \citep{parpas2015}, we take a different approach to reduce sample size.

Risk is an important component of optimization problems in many fields.
Mean-risk models are well-studied bi-objective models that allow decision makers to balance mean performance and risk aversion of some quantity (see, e.g., \citet{markowitz1959}, \citet{zenios1993}, \citet{ahmed2006}). 
Risk measures are tools to quantify risk and can be incorporated into stochastic programs to model risk aversion; however, modeling challenges due to severe, tail-of-distribution events remain.
We reference Chapter 6 in \citep{shapiro_lectures_2009} as well as the survey \citep{krokhmal2011} for thorough presentations of methods for modeling risk in stochastic programs.
A unique feature of our model relative to standard mean-risk models in stochastic programming  is that the mean and risk objectives may optimize different quantities corresponding to system operation in ``normal'' and ``extreme'' events.  To the best of our knowledge, the first work to consider separate models for normal and extreme events is \citep{Liu2016}. Their model differs from ours in that they include a chance constraint that the ``extreme'' events should occur with limited probability but otherwise consider expected value of both objectives, and sampling is not considered.

As discussed, the presence of two competing objectives is a key feature of our approach.
Standard ways for solving bi-objective, and, more generally, multi-objective, problems involve solving many slightly different optimization problems to find different trade-offs between objectives. 
We reference \citet{marler2004} for a survey of such techniques. 
In a survey of methods for solving stochastic, multi-objective optimization problems, \citet{gutjahr2016} present approaches to finding solutions under different sets of assumptions and problem types.
While we use a standard approach (the $\epsilon$-constraint method) for finding solutions with different trade-offs, such approaches would be useful for solving more computationally challenging extensions of our model. To the best of our knowledge, our sufficient condition for when replacing one objective in a bi-objective problem with an approximation can lead to a well-defined approximation of the set of Pareto optimal solutions is a new contribution to the bi-objective optimization literature.

\subsection{Summary of contributions}
We summarize our contributions (and outline the remainder of the paper) as follows:
\begin{itemize}
    \item We propose a bi-objective stochastic model for capacity expansion considering temperature as a key uncertain parameter (Section \ref{sec:modeling-approach}),
    \item We present solution methods including a conditional sampling approach to enhance results achieved when using SAA to solve the proposed model (Section \ref{sec:solution-approach}),
    \item We propose a spatial temperature modeling approach that can efficiently generate conditionally extreme temperature scenarios (Section \ref{sec:temp-model}),
    \item We conduct an extensive case study using realistic data that demonstrates the value of the proposed bi-objective model and the potential importance of incorporating temperature effects when capacity planning (Sections \ref{sec:casestudy} and \ref{sec:analyses}).
\end{itemize}

\section{Models}\label{sec:modeling-approach}
We describe two approaches for modeling the problem of generation expansion considering the uncertainty in demands and
generator outages. The approaches are both two-stage models having the same set of first-stage decision variables used
to decide on the generator fleet, described in Section \ref{sec:gen-decisions}. The two approaches also use the same
model of uncertain outcomes, described in Section \ref{sec:uncertain}. The first approach, described in Section
\ref{base-model-section} and referred to as the base model, is a traditional two-stage stochastic programming model that
minimizes the expected cost of the system operation. The second approach is our proposed bi-objective model, which we
present in Section \ref{sec:biobj-model}.  In order to focus on the challenges of achieving solutions with a desirable trade-off between average performance and risk in extreme events, we use a simple model of power grid operations in these models. In Section \ref{sec:model-extensions} we discuss extensions that could be incorporated within this modeling framework to make the power grid model more realistic.

\subsection{Generator decisions}\label{sec:gen-decisions}
The main purpose of the models is to choose the generator fleet in the system.
We consider two types of candidate generators, indexed by sets $\gensbin$ and $\genscont$, with $\gens = \gensbin \cup
\genscont$ defined as the set of all candidate generators. Generators in the set $\gensbin$ represent potential larger,
conventional generators (e.g., nuclear, hydro, and
fossil plants) for which the decision is to either build this generator or not.   
Thus, for each potential generator $i \in \gensbin$, we introduce a binary variable $\sone_i$ with $\sone_i=1$ indicating that generator $i$ is built.
Generators in the set $\genscont$ represent candidate generators that have a flexible capacity size (e.g.,
wind and solar farms where a large number of units could be installed). For each potential generator $i \in \genscont$ we introduce a continuous decision variable
representing the amount of generation capacity installed. Each candidate generator has a known location, so that generators with the same specifications but at different potential locations
are treated as different generators in the set of candidates.

The fixed costs of the generator decisions (for both building and operating the generators) are captured by the vector
$\costcap$: for a generator $i \in \gensbin$, $\costcap_i$ is the cost to build generator $i$ and have it as part of the
fleet; for a generator $i \in \genscont$, $\costcap_i$ is a cost per unit of generation capacity installed. These costs
are amortized to a cost per period (i.e., cost per hour in the case study, discussed in Section \ref{sec:casestudy}) so
that they can be added with the average operating cost per hour. Thus, the total fixed cost of a fleet represented by
the vector of generator decisions $\sone$ is
$\costcap^\top x$. Each generator $i \in \gens$ also has an associated generation cost function $\costgen_i(\cdot)$ such
that $\costgen_i(\stwo)$ represents the cost to produce $\stwo$ units of power with generator $i$ during a time period.

We let $X \subseteq \{0,1\}^{\gensbin} \times \mathbb{R}_+^{\genscont}$ represent the set of all feasible generator fleets. In addition to the binary and nonnegativity constraints on the decision variables, $X$ includes other constraints such as the upper bounds on the size of the generators $i \in \genscont$.

\subsection{Uncertainty modeling}
\label{sec:uncertain}
The generator fleet decisions must be evaluated across the full range of possible scenarios they may encounter.
We thus consider  performance across a set of time periods throughout the entire year as well as across the distribution
of  values of uncertain outcomes, including temperature, demand, and generator outages. To model variation of system
impacts throughout the year, we select a finite set of time periods $\timeperiods$. In our case study, for instance, we
use selected hours throughout the year that cover different seasons and times of day. We model the uncertain outcomes of
temperature, demand, and generator outages as random variables. In order to mathematically treat the evaluation of
performance across the set of time periods in the same way as over outcomes of these random variables, we define a
random variable $\Time$ which is supported on the set $\timeperiods$ and is uniformly distributed over time periods $\time \in \timeperiods$.

The random variables representing temperature, demand, and generator outages depend on $\Time$. The dependence of
temperature on the hour $\Time$ is handled by constructing separate spatial distributions of temperature for each hour
$\time \in \timeperiods$. These distributions are represented by the random variable $\Temp$ and are assumed to capture
spatial correlations among temperatures in different locations. In the context of the power grid, the correlations
across space are important for understanding the overall system capacity, especially during extreme temperature events
since these tend to have a wide geographic impact. Section \ref{sec:temp-model} describes an approach to creating such a
statistical model.  Demand at each location is modeled using linear regression that considers hour of day (as determined
by $\Time$), season of year (as determined by $\Time$), and temperature at that location as independent variables; the total demand $\demand$ is then obtained as a sum of the demands across all locations and represented as the random variable $\Demand$. 

We assume that, given a time observation $\time$, each generator $i \in \gens$ has a baseline capacity given by the function $\cap_i(\time)$. For generators $i \in \gensbin$, $\cap_i\paren{\time}$ can be interpreted as the nameplate capacity and typically does not vary with $\time$. For generators $i \in \genscont$, $\cap_i\paren{\time}$ is a factor of total capacity typically available at time $\time$. For example, solar farm production is largely a function of the position of the sun, which is determined by $\time$ (both hours of the day and time of year play a role, and $\time$ determines both). Wind farms also have daily and seasonal production patterns that are functions of $\time$. 

Uncertainty in generator capacities is represented by the random variables $\Outages_i$ for $i \in \gens$. 
For any generator $i$, $\Outages_i$ captures the randomness in its available capacity that occurs outside its baseline value $p_i(H)$. 
For instance, the availability of a thermal generator $i$ is modeled by a Bernoulli random variable $\Outages_i$
representing whether or not the generator is available, and the probability that $\Outages_i=1$ is a function of the temperature at generator $i$. 
For a renewable generator $i$ (like a wind farm), $\Outages_i$ is a continuous random variable representing uncertainty in the capacity factor (e.g., due to temperature, cloud cover, wind speed, etc.).
Therefore, the capacity of any generator $i$ in a time period $\Time$ is given by $\sone_i \cap_i(\Time)\Outages_i$. A detailed example of this method will be described in Section \ref{sec:casestudy}.

To simplify notation in our models, we let $\Rand = (\Time, \Temp, \Outages, \Demand)$ and call its distribution $\dist$.

\subsection{Base model}\label{base-model-section}
 In the base model, the second-stage cost is represented by the function $\obja$, which evaluates the operational cost of a design vector $\sone$ for a given observation of the random parameters $\xi=(\time,\temp, \outages, \demand)$:
\begin{subequations}\label{eq:qnormal}
\begin{align}
    \obja (\sone,\xi) := \min_{\stwo,\shed} &\ \sum_{i\in \gens} \costgen_i (\stwo_i) + \costshed \shed, \label{eq:qnormal:objetive}\\
    \text{s.t. }& \sum_{i \in \gens} \stwo_i \geq \demand - \shed, \label{eq:qnormal:demand}\\
        & \stwo_i \leq \sone_i \outages_i \cap_i(\time),\quad \quad \forall i \in \gens, \label{eq:qnormal:generators}\\
        & \stwo,\shed \geq 0. \label{eq:qnormal:nonneg}
\end{align}
\end{subequations}
In \eqref{eq:qnormal}, the decision variables $\stwo$ and $\shed$ are the recourse decisions made in response to a given observation of the parameters $\xi=(\time,\temp,\outages,\demand)$. The variables $\stwo_i$ for a generator $i \in \gens$ represent the power level of generator $i$, and the variable $\shed$ represents the level of load shed (unmet demand). The parameter $\costshed$ represents the penalty per unit of shed load, which sets the model's preference between cost and load shed.

The two-stage stochastic programming formulation is then as follows:
\begin{align}\label{eq:base-model}
    \min_{\sone \in X}  \costcap^\top \sone \ + & \ \E_\dist \bracket{\obja(\sone,\Rand)}.
\end{align} 
Here, we choose a generator fleet that minimizes the fixed cost of the fleet (amortized to a per-period basis) plus the expected second-stage cost, where the expectation captures the average across different time periods throughout the year as well as over the uncertain temperatures, generator availability, and demands. 

As discussed in the introduction, a modeling goal of this problem is to find multiple solutions with different trade-offs between cost and risk. By solving \eqref{eq:base-model} with various values of $\costshed$ in \eqref{eq:qnormal:objetive}, the modeler can change this risk preference to find solutions that weigh cost and load shed differently.

In the context of rare, high-impact events (e.g., a wave of extreme cold temperatures across an entire region that simultaneously reduces generator capacity and increases total demand, leading to significant unmet demand), this formulation has potential challenges. For one, approximating the expectation $\E_\dist \bracket{\obja(\sone,\Rand)}$ with some rare, high-impact $\rand$  via Monte Carlo methods would require a very large sample in order to even observe those rare events, making the task of optimizing this value very challenging. 
Another limitation of this formulation is that the second-stage recourse problem may actually change in extreme
situations. For example, load shed, which is deemed unacceptable in normal conditions, may be unavoidable under extreme
circumstances, and in such situations, the objective may change from minimizing cost to minimizing load shed.
Furthermore, due to the severe impact that large-scale load shedding can have, system operators may find that minimizing the
average cost leads to solutions in which the risk of very bad events is unacceptably high. 

\subsection{Bi-objective model}\label{sec:biobj-model}
Increasing the investment in generation capacity can help reduce risk from extreme events but at the same time
increases investment cost. This makes low cost and low risk from extreme events competing objectives. To make this reality explicit, we propose a
bi-objective model where the first objective is to minimize the expected total cost (i.e., investment plus operating
costs) of ``normal'' operation, and the second objective is to minimize a measure of the risk of ``extreme'' impacts. A
key benefit of this modeling framework is that it does not require  estimating the \textit{cost} of high-impact events like the base model \eqref{eq:base-model} does. While the first objective is understood as a cost, the second objective can have different units that need not be understood as a cost, as we demonstrate in the case study. We use the function $\obja$ to assess the cost of normal operation, so our first objective is \eqref{eq:base-model}.

To model the second objective of minimizing a measure of risk in extreme events, we define the problem $\objb$ to assess the impact of extreme events by estimating the minimum amount of load shedding required given the fleet $\sone$ and the scenario $\xi=(\time, \temp, \outages, \demand)$:
\begin{subequations}\label{eq:qextremelp}
\begin{align}
    \objb(\sone, \xi) :=\min_{\stwo,\shed} & \  \shed, \label{eq:qextremelp:objective}\\
    \text{s.t.  } & \sum_{i \in \gens} \stwo_i \geq \demand - \shed, \label{eq:qextremelp:demand}\\
        & \stwo_i \leq \sone_i \outages_i \cap_i(\time), \quad \quad \forall i \in \gens, \label{eq:qextremelp:generators} \\
        & \stwo, \shed \geq 0. \label{eq:qextremelp:nonneg}
\end{align}
\end{subequations}
The only difference between model \eqref{eq:qnormal} for normal operation and model \eqref{eq:qextremelp} for extreme events is the objective function: for extreme events, we care only about meeting load (with no regard for cost). As a result, we will always choose to dispatch all available generation capacity to cover the electric load, and we can simplify \eqref{eq:qextremelp} as follows:
\begin{equation}\label{eq:qextrememax}
\objb(\sone,\rand) = \paren{D - \sum_{i \in \gens}\sone_i\outages_i \cap_i(\time)}_+,
\end{equation}
where $(\cdot)_+ = \max\{\cdot,0\}$. 

To model risk aversion, we use a risk measure $\rho$ on $\objb$, so our second objective is
\begin{equation}\label{eq:objtwo}
    \min_{\sone \in X}\ \rho_\dist \bracket{\objb (\sone,\Rand)},
\end{equation}
where the notation $\rho_\dist$ is used to indicate that the risk measure $\rho$ is evaluated using the distribution $\dist$ of the uncertain outcomes.

Since the primary goal of the power system is to meet customer demand, power grid reliability is typically assessed by some combination of how much and how frequently demand is not met \citep{epri-ra}. 
The probabilistic risk measure loss-of-load expectation (LOLE) and its relative loss-of-load probability (LOLP) have been the incumbent metrics for assessing power grid reliability for decades \citep{billinton2015}. 
Other risk measures, such as conditional value at risk ($\cvar{}$), can assess power grid performance by considering the amount of load shed that occurs in the worst cases. 
Our definition of $\objb$ accommodates the use of these risk measures (among others) as $\rho$, and we specifically consider two of them.

The first risk measure is the conditional value at risk with threshold $\cvarthres$ ($\cvar[\cvarthres]$) of $\objb$. This is defined as
\begin{equation}\label{eq:cvar-formally}
    \cvar[\cvarthres]_\dist [\objb (\sone,\Rand)] = \inf_\gamma \braces{\gamma + \frac{1}{1-\cvarthres} \E \bracket{\paren{\objb\paren{\sone,\Rand} - \gamma}_+}},
\end{equation}
and in our setting it can be interpreted as the expected value of load shed in the worst $\cvarthres$ percentage of outcomes \citep{rockafellar2000}.

The second risk measure is the loss of load probability (LOLP), defined as the probability of load shed occurring. Given our definition of $\objb$, we can write this as
\begin{equation}\label{eq:lolp-formally}
    \lolp_\dist\bracket{\objb (\sone,\Rand)} = \P\bracket{ \objb (\sone,\Rand) > 0}.
\end{equation}

Incorporating $\cvar{}$ into an optimization model can be formulated with linear constraints and is theoretically supported as a coherent risk measure \citep{rockafellar2000}. See \citet{rockafellar2000} and \citet{ruszczynski2006} for more details on $\cvar{}$ and coherent risk measures. Formulating \lolp\ requires additional binary variables, so it is not convex and is more computationally challenging than $\cvar{}$. Therefore, in Section \ref{sec:analyses} we use $\cvar{}$ in most case study analyses, and we also directly compare $\cvar{}$ and \lolp.

Regardless of the choice of risk measure, the bi-objective model is
\begin{equation}\label{eq:biobj-abstract}
    \min_{\sone \in X}\ \paren{\costcap^\top \sone \ + \ \E_\dist \bracket{\obja(\sone,\Rand)},\ \rho_\dist \bracket{\objb (\sone,\Rand)}}.
\end{equation}
We note this is different from a common use of a risk measure in stochastic programming because the subproblem in the risk measure ($\objb$) is different from the subproblem in the expected value ($\obja$). Our implementation of this method in the case study is discussed in detail in Section \ref{sec:casestudy}. 

Considering that we have two objectives to optimize in \eqref{eq:biobj-abstract}, the appropriate solution concept is to (approximately) identify the set of \textit{Pareto optimal solutions} \citep{marler2004}. 

\begin{definition}
\label{def:pareto}
    A solution $\vone \in X$ is \textbf{Pareto optimal} to the bi-objective problem of minimizing the two functions $f$ and $g$ over the set of feasible solutions $X$ if there does not exist a solution $\vtwo \in X$ such that $f(\vtwo) \leq f(\vone)$, $g(\vtwo) \leq g(\vone)$, and $f(\vtwo) < f(\vone)$ or $g(\vtwo) < g(\vone).$
\end{definition}
In other words, a solution is Pareto optimal if there is not another solution that improves one objective without detriment to another objective.
Thus, the set of all Pareto optimal solutions shows the set of feasible non-dominated trade-offs between the two
objectives. Given access to this information, the decision-maker can then choose among the non-dominated solutions. We discuss approaches to finding Pareto optimal solutions in Section \ref{sec:solution-approach}.

A consequence of the bi-objective approach is that the purpose of $\costshed$ in $\obja$ changes. In the base model
\eqref{eq:base-model}, $\costshed$ is a parameter that is varied to find solutions with different trade-offs between cost and risk. In the bi-objective model \eqref{eq:biobj-abstract}, such solutions are the Pareto optimal solutions which can be found without varying $\costshed$. Therefore, $\costshed$ in \eqref{eq:biobj-abstract} is a fixed parameter that can be interpreted as the cost of unmet demand or the cost of buying power from the external market, and it should be estimated accordingly.

\subsection{Notation summary}
Table \ref{tab:notation} contains a summary of the notation used in our models.
\begin{table}[h]
    \begin{tabularx}{\linewidth}{ c  X }
    Symbol & Definition\\ \hline
    $\gensbin$ & set of generators with binary decisions\\
    $\genscont$ & set of generators with continuous decisions\\
    $\gens$ &  $\gensbin \cup \genscont$, set of all generators\\
    $\timeperiods$ & set of time periods throughout the year\\
    $\sone$ & vector of first-stage decision variables for generators (construction)\\
    $\stwo$ & vector of second-stage decision variables for generators (power level)\\
    $\shed$ & second-stage decision variable for amount of load shed\\
    $\Time$ & uniformly distributed random variable supported on set of time periods $\timeperiods$ \\
    $\Temp$ & spatially correlated, time-dependent random vector of location-specific temperatures\\
    $\Outages$ & vector of random capacity variations for generators\\
    $\Demand$ & random electricity demand \\
    $\Rand$ & short-hand for $\paren{\Time, \Temp, \Outages, \Demand}$ \\
    $\dist$ & distribution of $\Rand$; i.e. $ \Rand \sim \dist$\\
    $\costcap$ & vector of generator fixed costs \\
    $\costgen_i(\cdot)$ & generator operating cost function for generator $i \in \gens$\\
    $\cap_i(\time)$ & baseline generator capacity at time $\time$ for generator $i \in \gens$\\
    $\costshed$ & penalty per unit of load shed\\ 

    \end{tabularx}
    \caption{Summary of notation.}
    \label{tab:notation}
\end{table}

\subsection{Model extensions}\label{sec:model-extensions}

The subproblems described in Sections \ref{base-model-section} and \ref{sec:biobj-model} used to evaluate the performance of a candidate generation fleet $x$
only consider the available capacity of generators and ignore many realistic aspects of the power grid. While we will
focus on those models in our case study, we describe in this section how this approach can be extended to accommodate a
variety of modifications to make the model more realistic. We leave the
implementation of these extensions to future work due to their computational challenge and our desire to focus the
current work on the value of the proposed bi-objective model and integration of a spatial temperature distribution into
the model.

\paragraph{More general operational models}
First, we can consider more general formulations of the operational subproblem $\obja$:
\begin{subequations}\label{eq:qnormal-general}
\begin{align}
    \obja(\sone,\rand) = \obja(\sone,\time,\temp,\outages,\demand) := \min_{\stwo,\recvara}\ & \costgena(\stwo) + \shedfn \paren{\recvara} , \\
    \text{s.t. } &0 \leq \stwo_i \leq x_i \outages_i \cap_i(\time) \quad \forall i \in \gens , \\
    & (\stwo,\recvara) \in \recfeasa(\time,\temp,\outages,\demand). \label{eq:qnormal-general-feasible-region}
\end{align}
\end{subequations}
Here, the recourse variables $\stwo_i$ for $i \in \gens$ are the generator power levels (as before), and the variable $\recvara$ is a vector of other recourse variables that could include load shedding, unit commitment decisions, transmission line power flows and topology, and other operational decisions desired by the modeler. The function $\costgena$ represents the cost to operate the generators $\gens$ at levels $\stwo$ while the function $\shedfn$ represents costs associated with the other recourse variables like transmission costs or short-term preventative measures. The set $\recfeasa$ represents other constraints associated with grid operations, corresponding to the variables $\recvara$, and could include constraints like AC or DC power flow equations, transmission line limits, and reserve requirements.

A model $\objb$ with the same structure as \eqref{eq:qnormal-general} can also be used for extreme events, though we emphasize that the additional recourse decision variables $\recvarb$, the objective function, and the feasible region $\recfeasb(\time,\temp,\outages,\demand)$ may all be different in this problem than in \eqref{eq:qnormal-general}, representing the reality that grid operators may change their operating constraints and goals when the system is facing extreme conditions. Possible modifications could include preemptively de-energizing power lines to mitigate wildfire risk and implementing short-term hardening measures in preparation for some event (e.g., putting out sandbags prior to a hurricane to mitigate flood risk). Such a definition of $\objb$ can be used with the $\cvar{\alpha}$ risk measure as before or with a slightly generalized version of the \lolp\ risk measure, denoted by $\chance$ with the parameter $\chanceparam$, which we define as the probability that $\objb$ exceeds some threshold $\chanceparam$. This can be written as
\begin{equation}\label{eq:chance-formally}
    \chance^\chanceparam_\dist\bracket{\objb (\sone,\Rand)} = \P\bracket{ \objb (\sone,\Rand) > \chanceparam}.
\end{equation}
With $\chanceparam = 0$ and certain definitions of $\objb$, this is equivalent to \lolp.

\paragraph{Transmission expansion}
While omitted from our case study, transmission network modeling is an important consideration in many electricity planning problems. In addition to modeling an existing network as mentioned above (to capture power flows and line capacities), the model can include transmission expansion decisions by including additional binary variables in the first stage variables $\sone$, representing potential new transmission lines whose parameters are included in \eqref{eq:qnormal-general-feasible-region}. 

\paragraph{Scenarios of different duration}
As written so far, time is discretized into standalone, hour-long periods represented by the random variable $\Time$.
Specifically, while this model considers a variation of time periods throughout the year, it ignores the
interrelationship between decisions made in consecutive time periods and the correlation of uncertain quantities over time. In an extension of the model, a single ``period'' $\time$  could instead represent a \textit{sequence} of time periods (minutes, half hours, etc.) connected through time so that decisions
made in a previous time period impact a later time period. In this extension, the evaluation problems \eqref{eq:qnormal} and \eqref{eq:qextremelp} would become multi-period optimization problems. The notation presented so far is intended to accomplish such
a definition of scenarios by extending to vector-valued functions. This extension would allow modeling more generator
operational constraints, like ramping rates, energy storage, and unit commitment decisions, all within
\eqref{eq:qnormal-general-feasible-region}. Thus, rather than using a selection of hours $H$ throughout the year, an extension of the model could, for example, use a selection of days throughout the year, where the second-stage problem associated with each day would include multiple time periods throughout the day that would allow capturing temporal dependence in temperatures and temporal constraints such as generator ramping limits.

\section{Solution Approach}\label{sec:solution-approach}
A key challenge in solving either the base model \eqref{eq:base-model} or the bi-objective model \eqref{eq:biobj-abstract} is the presence of the expected value, which appears in both problems, and the risk measure, which appears only in the bi-objective model.
Since evaluating such terms exactly is, in general, intractable, we focus on methods that approximate these with a set of scenarios that may be obtained by sampling observations from $\Rand$. We now discuss the details of this for both the base model and the bi-objective model.

\subsection{Base model}
\label{sec:basesolve}
To solve the base model \eqref{eq:base-model}, we can use SAA to approximate the expected value \citep{saa2002}. Specifically, a random sample of observations of the uncertain data $\{\rand^k, k \in \N\}$ is taken and the expected value in \eqref{eq:base-model} is replaced by an average over the observations in this sample. Thus, the sample average approximation of the base model becomes
\begin{align}\label{eq:base-model-saa}
    \min_{\sone \in X} \ &  \costcap^\top \sone + \frac{1}{|\N|}\sum_{k\in \N} \obja \paren{\sone,\rand^k}.
\end{align}

It is important to recognize that the solution obtained from \eqref{eq:base-model-saa} is random due to its dependence on the random sample. Thus, a standard approach in applying SAA is to solve multiple replications to obtain multiple candidate solutions. In addition, the optimal value of the SAA problem is a biased estimate of the objective value of the solution, since the solution is obtained by optimizing to that sample. Thus, after obtaining a solution, it is necessary to evaluate the solution on a new, independent, sample to obtain an unbiased estimate of its cost. Finally, recall that our goal is to find solutions that explore the trade-off between the risk and cost objectives. 

With these considerations in mind, we propose a two-phase approach for using the base model to obtain solutions that approximate the efficient frontier between the risk and cost objectives:
\begin{enumerate}
\item \textit{Solution generation}: solve \eqref{eq:base-model-saa} over multiple independent samples and many values of the parameter $\costshed$ (the cost per unit of unmet demand) to generate multiple first stage solutions having different trade-offs between average cost and average load shed.
\item \textit{Solution evaluation}: use an independent sample to estimate the cost objective \eqref{eq:base-model} and the risk objective \eqref{eq:objtwo} for each of the obtained solutions. 
\end{enumerate}
The non-dominated pairs of estimated objective values from the second phase are then used as an approximation of the trade-off curve between the two objectives. We emphasize that the sample size in the solution evaluation phase may be taken to be much larger than that in the solution generation phase, since in that case we are evaluating the two objectives on the fixed set of generated solutions, as opposed to attempting to find an optimal solution to the SAA \eqref{eq:base-model-saa}.

Considering that it will be difficult to solve \eqref{eq:base-model-saa} with a huge number of scenarios, a potential limitation of \eqref{eq:base-model-saa} is that it may fail to find solutions that adequately protect against the effects of rare, high-impact events, since without a huge sample size such events may not even be observed.
Additionally, results from the literature studying the
convergence of SAA provide estimates on the sample size required to obtain reasonable approximations
\citep{shapiro2009}. These estimates demonstrate that the required sample size grows as the variance of the term
appearing in the expectation being approximated grows. Given that rare, high-impact events occur in $\Rand$, the
objective in \eqref{eq:base-model-saa} will have high variance. As a result, large samples of $\Rand$ may also be
required to obtain good results using SAA. We propose an approach to alleviate this sample size issue in Section \ref{sec:biobj-conditional}.

There are a variety of methods for solving the SAA problems \eqref{eq:base-model-saa}. In our computational study, we solve it by building the extensive form \citep{shapiro2009} and solving it with a general-purpose MIP solver Gurobi. For completeness, we describe the extensive form in the Appendix \ref{app:forms}. Other methods, such as Benders decomposition \citep{benders196,vanslyke1969} could be used instead. 

\subsection{Bi-objective: Unconditional}\label{sec:biobj-unconditional}
For the bi-objective model \eqref{eq:biobj-abstract}, we can also apply SAA by creating two samples from the same
distribution: $\None$ is used for estimating the first objective and $\Ntwo$ for estimating the second objective. (We
could use the same sample for both objectives, but we might want a larger sample for the second objective since it may
be more difficult to estimate.) Since the first objective is the same as the objective in the base model, the approach will be identical: we will replace the expected value with a sum over $|\None|$ observations of $\Rand$. The second objective is a little different from the first due to the risk measure, but we can still use a similar approach adapted to the risk measure we choose. To use $\cvar{}$, we first define $\cvar{}$ for a sample $\Ntwo$ from $\Rand$:
\begin{equation}\label{eq:cvar-discrete}
    \cvar[\cvarthres] \bracket{\objb \paren{\sone,\Rand^{\Ntwo}}} = \inf_\gamma \braces{\gamma + \frac{1}{1-\cvarthres} \sum_{k \in \Ntwo} \bracket{\paren{\objb\paren{\sone,\rand^k} - \gamma}_+}},
\end{equation}
Then the sample average approximation of the bi-objective model is 
\begin{equation} \label{eq:biobj-unc-cvar}
    \min_{\sone \in X} \ \paren{
    \costcap^\top \sone + \frac{1}{|\None|}\sum_{k\in \None} \obja \paren{\sone,\rand^k}, \
    \cvar[\cvarthres] \bracket{\objb \paren{\sone,\Rand^{\Ntwo}}}
    }.
\end{equation}

As discussed in Section \ref{sec:biobj-model}, solving this bi-objective problem \eqref{eq:biobj-unc-cvar} means finding
the set of Pareto optimal solutions. A common way to find Pareto solutions is the weighted sum method, where the objective function is the weighted sum of the two objectives \citep{das1997}. By solving the problem for different weightings of objectives, different Pareto solutions can be found. 
A drawback of the weighted sum method is that, when applied to a problem with a nonconvex feasible region, it does not find any Pareto solutions in a nonconvex part of the Pareto curve \citep{das1997}. Since our feasible region is nonconvex due to binary decision variables, we apply the $\epsilon$-constraint method to generate Pareto solutions of our bi-objective model \citep{marler2004}. 
Specifically, for varying values $U$ representing an upper bound on the second objective, we solve the following problem:
\begin{subequations}\label{eq:biobj-unc-cvar-epsilon}
\begin{align}
    \min_{\sone \in X} \ &  \costcap^\top \sone + \frac{1}{|\None|}\sum_{k\in \None} \obja \paren{\sone,\rand^k} \label{eq:biobj-unc-cvar-epsilon:objective}\\
    \st & \cvar[\cvarthres] \bracket{\objb \paren{\sone,\Rand^{\Ntwo}}} \leq \bound. \label{eq:biobj-unc-cvar-epsilon:constraint}
\end{align}
\end{subequations}
Solving this problem with different values of $\bound$ will find solutions with different trade-offs between the two objectives.

To approximate the Pareto frontier, we propose to use the same two-phase approach as described in Section \ref{sec:basesolve} for the base model.
For phase one, we solve \eqref{eq:biobj-unc-cvar-epsilon} over multiple different samples and  values of $\bound$ to generate a collection of solutions.
Then, we complete phase two exactly as described in Section \ref{sec:basesolve} to obtain unbiased estimates of the two objectives for each solution and use the solutions with non-dominated estimates as the approximation of the set of Pareto optimal solutions.
While we use the $\epsilon$-constraint method for simplicity in our case study, other methods for generating Pareto solutions (phase one) exist that are able to better approximate the set of Pareto solutions with less computational effort (see \citet{halffmann2022} and \citet{marler2010} for more details).

We can also use \lolp\ as the risk measure instead of $\cvar{}$ with the following bi-objective model:
\begin{equation}\label{eq:biobj-unc-lolp}
    \min_{\sone \in X} \ \paren{
    \costcap^\top \sone + \frac{1}{|\None|}\sum_{k\in \None} \obja \paren{\sone,\rand^k}, \
    \frac{1}{|\Ntwo|} \sum_{k\in \Ntwo} \indicator \bracket{\objb\paren{\sone,\rand^k} > 0} },
\end{equation}
where $\indicator[\cdot]$ is an indicator function defined to be 1 if the input is true and 0 otherwise. We can apply the $\epsilon$-constraint method in a similar way for this model as well.

With either risk measure on the second objective, we expect to find solutions that do better in the extremes. However, the sample size challenge of solving \eqref{eq:base-model-saa} remains since we still require large samples to see rare, high-impact scenarios.

\subsection{Bi-objective: Conditional} \label{sec:biobj-conditional}
To address the potential sample size challenge of solving \eqref{eq:biobj-unc-cvar} or \eqref{eq:biobj-unc-lolp}, we propose a \textit{conditional} model in which we approximate the risk objective with scenarios generated from a suitably chosen extreme random variable  $\RandE \sim \distE$, where $\distE$ is a conditional distribution of $\Rand$ given that $\Rand$ is ``extreme'' by some suitable definition, and $\rho'$ is a risk measure that may be different than $\rho$. Then we propose using the following model as an approximation of \eqref{eq:biobj-abstract}:
\begin{equation}\label{eq:biobj-cond-abstract}
    \min_{\sone \in X}\ \paren{\costcap^\top \sone \ + \ \E_\dist \bracket{\obja(\sone,\Rand)},\ \rho'_{\distE} \bracket{\objb (\sone,\RandE)}}.
\end{equation}
The motivation for this proposal is that with relatively few scenarios, a sample-average approximation of the second objective may capture the effects of rare, high-impact events.

The effectiveness of this idea depends on the definition of ``extreme'', as determined by the choice of $\distE$, and the choice of $\rho'$ such that
$\rho'_{\distE}\bracket{\objb\paren{\sone,\RandE}}$ is an adequate approximation of
$\rho_\dist\bracket{\objb\paren{\sone,\Rand}}$. We advocate choosing $\distE$ based on problem insight, specifically,
based on knowledge of the types of scenarios that usually lead to the worst outcomes for the objective $\objb
\paren{\sone,\rand}$. In the power grid application, this means conditioning on events that tend to lead to the most
significant power shortfalls. In our model, extreme temperatures, both high and low, can lead to a combination of high
power demand and reduced capacity, and hence tend to be the outcomes that lead to high power shortfalls. Thus, in our
case study we define the extreme distribution to be conditional that the spatial average temperature of the entire
region is in the hottest or coldest 1\% of the distribution of such spatial average temperatures. While we use this definition, our framework can accommodate alternative definitions. For example, one could define ``extreme'' as being events that are known to lead to bad outcomes on some reference solution, such as the existing system.

Since this extreme distribution $\distE$ is expected to contain a significantly higher fraction of bad
outcomes for $\objb \paren{\sone,\randE}$, we can choose $\rho'$ to be much less focused on the extreme worst cases than
$\rho$. For example, when using $\cvar{\cvarthres}$ as the risk measure $\rho$ in \eqref{eq:biobj-abstract}, we use $\cvar{\cvarthres'}$ as $\rho'$ in \eqref{eq:biobj-cond-abstract} where $\cvarthres' \gg \cvarthres$. A major benefit of this change is that when using SAA, approximating the risk measure with $\cvarthres'$ requires fewer scenarios than the risk measure with $\cvarthres$. Determining an appropriate value of $\cvarthres'$ may be challenging, but our results indicate that the approach is not sensitive to the particular choice of $\cvarthres'$. 

A key aspect of this approach that makes it computationally feasible is that the extreme distribution $\distE$ is chosen \textit{a priori}, before solving the optimization model. This enables the use of SAA to approximately solve \eqref{eq:biobj-cond-abstract} in the same way we solve \eqref{eq:biobj-abstract} by using samples of $\RandE$ (denoted $\Ntwo^E$) instead of $\Rand$ for the second objective. 
With $\cvar{}$ as the risk measure, this means we solve the following problem, which is an approximation of $\eqref{eq:biobj-unc-cvar}:$
\begin{equation} \label{eq:biobj-cond-cvar}
    \min_{\sone \in X} \ \paren{
    \costcap^\top \sone + \frac{1}{|\None|}\sum_{k\in \None} \obja \paren{\sone,\rand^k}, \
    \cvar{{\cvarthres'}} \bracket{\objb \paren{\sone,\RandE^{\Ntwo^E}}}
    }.
\end{equation}

If \lolp\ is used as the risk measure, the conditional sampling approach solves the following model to approximate \eqref{eq:biobj-unc-lolp}:
\begin{equation}\label{eq:biobj-cond-lolp}
    \min_{\sone \in X} \ \paren{
    \costcap^\top \sone + \frac{1}{|\None|}\sum_{k\in \None} \obja \paren{\sone,\rand^k}, \
    \frac{1}{|\Ntwo^E|} \sum_{k\in \Ntwo^E} \indicator \bracket{\objb\paren{\sone,\randE^k} > 0} }.
\end{equation}
While there is no analogous parameter $\cvarthres$ to adjust for \lolp, we can allow larger \lolp\ values in \eqref{eq:biobj-cond-lolp} than we would in \eqref{eq:biobj-unc-lolp} to achieve the same level of risk of extreme events.

We once again employ the two-phase method described in Section \ref{sec:basesolve} to use such sample approximations to approximate the efficient frontier.
In phase one, we generate solutions to the conditional model (either \eqref{eq:biobj-cond-cvar} or \eqref{eq:biobj-cond-lolp}), again using the $\epsilon$-constraint method to generate solutions with different trade-offs between the two objectives and solving each with multiple samples to generate many candidate solutions. 
In phase two, we evaluate all the candidate solutions exactly as described Section \ref{sec:basesolve}, where we emphasize that in the evaluation phase we estimate the \textit{true} risk measure $\rho$ using samples of the unconditional distribution $\Rand$. As discussed in the last section, it is generally feasible to use a suitably large sample for the evaluation phase because this is done for just a fixed set of solutions, rather than for the purpose of finding a solution.

The requirement that the extreme distribution be chosen before solving the optimization model means that the choice of extreme distribution cannot depend on decisions the optimization model makes. Therefore, the distribution will inevitably not completely match the scenarios where the most extreme outcomes occur for every candidate fleet. For example, if a particular candidate fleet has a lot of solar capacity, the risk is probably lower during hours when the sun is up compared to the risk for a comparably-sized fleet with little solar. While necessary for the approach, this does introduce error between the bi-objective model \eqref{eq:biobj-abstract} and the approximation \eqref{eq:biobj-unc-cvar}.

\subsubsection*{Approximation of Pareto solutions.} The general idea of this conditional sampling approach is to replace one objective in the bi-objective model that is difficult to optimize (the risk measure in extreme events) with an alternative objective that in some way approximates this but is easier to handle (the risk measure evaluated on the conditional distribution). We next explore when such an approximation scheme for a bi-objective problem may be expected to work well in the sense that it yields a set of solutions that approximate the true set of Pareto optimal solutions.

Abstractly, we consider a bi-objective problem with two objectives, $f$ and $g$ defined over some feasible region $X$:
\begin{align}\label{eq:min-f-g} 
\min_{x \in X} \ & \paren{f(x),g(x)}, 
\end{align}
and we define $X^*$ to be the set of Pareto optimal solutions of \eqref{eq:min-f-g}.  Our goal is to find $X^*$, but we assume solving this problem is computationally challenging due to the function $g$.
 In our model, we can think of $f(x)$ as 
$\costcap^\top \sone + \E_\dist \bracket{ \obja \paren{\sone,\Rand}}$ and $g(x)$ as $\rho_\dist \bracket{\objb (\sone,\Rand)}$.
We assume we are able to compute $\hat{X}$, which is the set of Pareto optimal solutions of the bi-objective problem:
\begin{align}
\min_{x \in X} \ & \paren{f(x),\hat{g}(x)} \label{eq:min-f-hatg} 
\end{align}
for a suitably chosen function $\hat{g}$, which in our context would be $\rho_{\distE}' \bracket{\objb \paren{\sone,\RandE}}$ (or possibly even a sample-average approximation of this term).

We next provide a relaxed version of the definition of Pareto optimal solutions given in Definition \ref{def:pareto}.
\begin{definition}
\label{def:epspareto}
    A solution $\vone \in X$ is $\epsilon$-\textbf{Pareto optimal} to the bi-objective problem \eqref{eq:min-f-g} if there does not exist a solution $\vtwo \in X$ such that $f(\vtwo) \leq f(\vone)$ and $g(\vtwo) < g(\vone)-\epsilon$.
\end{definition}

Intuitively, if $\hat{g}(x) \approx g(x)$ for all $x \in X$ we would expect the Pareto optimal solutions of \eqref{eq:min-f-hatg} to be a good approximation of the set of Pareto optimal solutions of \eqref{eq:min-f-g}. The following proposition shows that a much weaker condition is sufficient: if $\hat{g}$ is sufficiently close to \textit{a monotone transformation of $g$}, then the set of Pareto optimal solutions $\hat{X}$ of \eqref{eq:min-f-hatg} is an approximation of the set of Pareto optimal solutions $X^*$ of \eqref{eq:min-f-g}. Specifically, under this assumption, every solution in $\hat{X}$ is approximately Pareto optimal to \eqref{eq:min-f-g}, and for every Pareto optimal solution to \eqref{eq:min-f-g} there is a solution in $\hat{X}$ that is nearly as good or better on both objectives.

\begin{proposition}\label{proposition:pareto-monotone}
Let $X^*$ be the set of Pareto-optimal solutions to \eqref{eq:min-f-g} and $\hat{X}$ be the set of Pareto-optimal solutions to \eqref{eq:min-f-hatg} and assume both these sets are non-empty.
Assume there exist a positive constant $L$, a differentiable monotone increasing function $h$ with $h'(y) \geq L $ for all $y \in [\min_{x \in X} g(x), \max_{x \in X} g(x)]$, and  $\epsilon > 0$ such that 
    \begin{equation}\left| \ghat(x) - h(g(x))\right| < \epsilon \quad \text{ for all } x \in X. \end{equation}
Then,
\begin{enumerate}
\item[(a)] Every $\hat{x} \in \hat{X}$ is a $2\epsilon/L$-Pareto optimal solution of \eqref{eq:min-f-g}.
\item[(b)] For every $x^* \in X^*$, there exists $\hat{x} \in \hat{X}$ which satisfies:
\begin{equation}
\label{eq:neardom}
f(\hat{x}) \leq f(x^*), \quad g(\hat{x}) \leq g(x^*) + 2\epsilon/L .
\end{equation}
\end{enumerate}
\end{proposition}

\begin{proof}
(a) Let $x \in \hat{X}$, and let $\bar{x} \in X$ be such that $f(\bar{x}) \leq f(x)$.
 Then $\ghat(x) \leq \ghat(\bar{x})$ since $x$ is Pareto optimal to \eqref{eq:min-f-hatg}. Thus,
    \[ h(g(x)) - h(g(\bar{x})) \leq \ghat(x) + \epsilon - \ghat(\bar{x}) + \epsilon \leq 2 \epsilon.\]
    Therefore,
    \begin{align*}
        h^{-1}(h(g(x))) &\leq h^{-1}(h(g(\bar{x})) + 2 \epsilon) \implies\quad
        g(x) \leq g(\bar{x}) + 2\epsilon/L.
    \end{align*}
    That is, $g(\bar{x}) \geq g(x) - 2 \epsilon/L$ and hence $x$ is $2\epsilon/L$-Pareto optimal.

(b)
We first observe that for any feasible solution $x \in X$ there exists a solution $\hat{x} \in \hat{X}$ such that $\hat{g}(\hat{x}) \leq \hat{g}(x)$ and $f(\hat{x}) \leq f(x)$. Indeed, the claim is trivial if $x \in \hat{X}$, and on the other hand if $x \notin \hat{X}$ then by definition of Pareto-optimality there exists $\hat{x} \in X$ such that $\hat{g}(\hat{x}) \leq g(x)$, $f(\hat{x}) \leq f(x)$, and one of those inequalities is strict. Among all such solutions, let $\hat{x}$ be one that minimizes $f(\hat{x})+\hat{g}(\hat{x})$. We argue that $\hat{x} \in \hat{X}$. If not, we can repeat the above argument to establish existence of $x' \in X$  which satisfies $f(x') \leq f(\hat{x}) (\leq f(x))$, $\hat{g}(x') \leq \hat{g}(\hat{x}) (\leq \hat{g}(x))$ and $f(x') + \hat{g}(x') < f(\hat{x}) + \hat{g}(\hat{x})$, which contradicts the choice of $\hat{x}$.

Let $x^* \in X^*$. By the observation of the previous paragraph, there is an $\hat{x} \in \hat{X}$ that satisfies:
\[ f(\hat{x}) \leq f(x^*), \quad \hat{g}(\xhat) \leq \hat{g}(x^*). \]
Then,
    \[ h(g(\xhat)) - h(g(x^*)) \leq \ghat(\xhat) + \epsilon - \ghat(x^*) + \epsilon \leq 2 \epsilon.\]
    Therefore,
    \begin{align*}
        h^{-1}(h(g(\xhat))) &\leq h^{-1}(h(g(x^*)) + 2 \epsilon) \implies\quad
        g(\xhat) \leq g(x^*) + 2\epsilon/L.
    \end{align*}
\end{proof}

Applying this result to our model, this implies that if we can choose $\rho_{\distE}' \bracket{\objb \paren{\sone,\RandE}}$ such that, \textit{after a monotone transform}, it approximates $\rho_\dist \bracket{\objb (\sone,\Rand)}$, then we may expect the set of Pareto solutions we obtain to perform well. In particular, even if $\rho_{\distE}' \bracket{\objb \paren{\sone,\RandE}}$ is significantly shifted or scaled from $\rho_\dist \bracket{\objb (\sone,\Rand)}$, it may still serve as a suitable proxy for the purpose of generating an approximation of the set of Pareto optimal solutions. In our numerical experiments in Section \ref{sec:conditionalapprox}, we demonstrate empirically that our choice of conditional distribution exhibits this property. 

We also observe that satisfying Proposition \ref{proposition:pareto-monotone} is sufficient but not necessary for \eqref{eq:biobj-cond-abstract} to be a good model of \eqref{eq:biobj-abstract}. For instance, the relationship between the unconditional and conditional risk measures at far-from-optimal solutions is not important since those solutions are unlikely to be selected as an optimal solution of the approximate problem.




\section{Spatial Temperature Modeling}\label{sec:temp-model}
We are interested in the effects of temperature, especially extreme temperature, on power grid operations in a large geographical region over the course of the year. To model these relationships, we first devise a method to generate appropriate temperature samples across our region of interest. Since electric load and generator availability are temperature sensitive (this is discussed in Section \ref{sec:temperature-dependence}), we desire samples that have a high enough spatial resolution to provide information at the geographical location of individual loads and generators, while also accounting for variations across the region and capturing the correlation of temperatures in locations that are close to each other. 

We also desire a time resolution of scenarios that can capture different behaviors at different times of the day, for two reasons. First, electricity supply and demand have well-defined patterns over the course of a day (summer afternoon demand peaks, overnight demand troughs, solar capacity factors, etc.). We also observe that it is the \emph{combination} of (i) extreme temperature, (ii) high electric load, and (iii) time-dependent generator availability (i.e., wind and solar), and not necessarily any one of them alone, that can lead to high load shed in the grid. Here, we explain how we create a statistical model and generate scenarios to use in the SAA and solution evaluation process described previously.


\subsection{Temperature Data}

The temperature dataset used to fit our statistical model is a reanalysis product from the North  American  Land  Data  Assimilation System project, phase 2 (NLDAS-2), with hourly values of two-meter above ground temperature in the Midwestern United States from 1991-2019 at a 1/8th degree latitude-longitude resolution \citep{xia2012continental1,xia2012continental2}.
Specifically, we consider data for the following states: Wisconsin, Iowa, Minnesota, Illinois, Indiana, and Michigan (including Upper Peninsula).
In total, there are 535 counties, and we define the temperature data for each county from the NLDAS product using the gridpoint closest to the centroid of each county.

\subsection{Statistical Model}\label{sec:stats-model}

In this section, we provide the statistical details of the spatial temperature model.
We use our model to simulate temperature fields from both ``normal'' and ``extreme'' scenarios, where ``normal'' is an unconditional simulation from the model, and ``extreme'' is a simulation conditional on the spatial average of all 535 county-level temperatures being above (below) its $0.99$ $(0.01)$ quantile. 
This formulation stands in contrast to common models in statistical extremes that only consider the upper tail of the distribution and ignore the lower tail and bulk of the distribution. A consequence of bulk-and-tails modeling is that including information from the bulk of the distribution can improve the accuracy of the fit in the tails and thus result in more accurate modeling of extremes \citep{stein2020}. 	
As demonstrated in Section \ref{sec:comparespatialvsind}, including spatial correlation in the simulated temperature fields is essential for obtaining solutions with lower risk in the power system investment model.

We employ a copula to model the spatial dependence in temperatures, thus allowing for separate modeling of the marginal distributions (at a single location) and joint distribution (across the region) of county-level temperatures.
A copula $C : [0,1]^n \to [0,1]$ is a cumulative distribution function (cdf) with uniform marginal distributions \citep{Nels06,Joe14}.
With $\Phi$ denoting the cdf of the standard Gaussian distribution and $\Phi_R$ the multivariate standard Gaussian cdf with correlation matrix $R$, the Gaussian copula equals 
\begin{equation}
\label{eq:gaussiancopula}
    C(u_1,\dots,u_n) = \Phi_R(\Phi^{-1}(u_1),\dots,\Phi^{-1}(u_n)).
\end{equation} 
By applying componentwise quantile functions to the Gaussian copula, which takes care of spatial dependence in the model, we produce a flexible multivariate distribution whose marginal distributions parametrically describe extreme temperatures (specifically, both hot and cold extremes, as well as the bulk of the temperature distribution).
We also devise an efficient sampler for spatial extremes which hinges upon the choice of Gaussian copula combined with our definition of extreme in terms of a spatial average exceeding a threshold.

For every hour $\time$ used in our model, we aim to produce a joint distribution of the temperatures across multiple prescribed spatial sites $l=1,\ldots,L$. Recall that the hour $\time$ refers to a particular hour in the year, which implies its hour within the day $h$ and the day of year $d \in \Dstat$, where $\Dstat$ is the set of days. Since our statistical model will need the hour-of-day $h$ and the day of year $d$, in this section we will specify the hour of year via $h$ and $d$. Thus, our aim is to estimate a joint distribution of the vector $\{T_{l,d,h}: 1 \leq l \leq L\}$, where $T_{l,d,h}$ is the temperature of location $l$ on day $d$ in hour-of-day $h$. 
The modeling is carried out separately for each hour. However, to reduce the computational burden, we only model the even hours of the day $h=2,4,\ldots, 24$ and perform the following steps separately for each hour $h$. Thus, in the remainder of this section we often suppress the index $h$. The model can easily be extended to all 24 hours by paying the corresponding computational cost. 

\paragraph{Modeling the marginal distributions}
\label{sec:marginalmodels}
For the marginal distribution of temperature at hour $h$ of day $d$ for site $l$, $1 \leq l \leq L$, we use a seven-parameter Bulk-And-Tails (BATs) distribution \citep{stein2020}. This model allows us to capture complex hot and cold tail behaviors, while simultaneously modeling the main range of temperatures (bulk). It can also account for seasonal dependence at each location using simple parameter embedding, which we describe below.  
For this setup, the cdf of the BATs marginal distribution of temperature $T_{l,d}$ at site $l$ on day $d$ at hour $h$ equals $F_{\nu_l}(\Hstat_{\theta_{l,d}}(T_{l,d}))$, where $F_{\nu_l}$ is the cdf of the student-$t$ distribution with $\nu_l$ degrees of freedom, and 
\[ \Hstat_{\theta_{l,d}}(T_{l,d})  = \left\{ 1 + \kappa_{1,l} \Psi \left( \frac{T_{l,d}-\phi_{1,l,d}}{\tau_{1,l,d}} \right) \right\}^{1/\kappa_{1,l}} - \left\{ 1 + \kappa_{0,l} \Psi \left( \frac{\phi_{0,l,d}-T_{l,d}}{\tau_{0,l,d}} \right) \right\}^{1/\kappa_{0,l}}.
\]
The parameters to be estimated are the smoothness parameter $\nu_l$ and $\theta_{l,d}=\{\phi_{0,l,d}$,\allowbreak$\phi_{1,l,d}$,\allowbreak$\tau_{0,l,d}$,\allowbreak$\tau_{1,l,d}$,\allowbreak$\kappa_{0,l}$,\allowbreak$\kappa_{1,l}\}$. Here $\phi_{0,l,d},\phi_{1,l,d}$ are location parameters, $\tau_{0,l,d}, \tau_{1,l,d}$ are nonnegative scale parameters, and $\kappa_{0,l},\kappa_{1,l}$ are shape parameters which control the tail behavior of the lower and upper tails of the distribution. As the subscripts indicate, the scale and location parameters depend on both day and location, whereas the shape parameter depends on location only; this choice is discussed in \citet{krock2022}.
That is, the location and scale parameters of the BATs distribution depend on seasonal covariates; the mappings $d \rightarrow \phi_{0,l,d},\phi_{1,l,d},\tau_{0,l,d}, \tau_{1,l,d}$ are represented by periodic cubic splines with nine coefficients, for a total of 36 to be fitted at each site $l$, $1 \leq l \leq L$. The model thus has $3L$ parameters for $\kappa_{0,l},\kappa_{1,l},\nu_l$, and $36L$ parameters for $\phi_{0,l,d},\phi_{1,l,d},\tau_{0,l,d}, \tau_{1,l,d}$ for each of the 12 target hours. This results in a total of $39L \times 12$ parameters to express all marginal distributions of temperature. However, this model is fitted separately per location and hour, resulting in maximum likelihood problems with 39 parameters carried out separately before any joint modeling is undertaken. The data used to fit this model is the daily NLDAS-2 data at the chosen hour $h$ in the 1991-2019 yearly range at the given site $l$, which results in $365 \times 29$ plus 7 leap years data points for a total of $10592$ points (the total number of days) to estimate 39 parameters.
See \citet{krock2022} for a similar analysis of daily temperature. 

\paragraph{Modeling the joint distribution through copula}

The county-level temperatures are marginally transformed to be standard normal random variables using the transformation 
\begin{equation}
\label{eq:definez}
    T_{l,d} \mapsto z_{l,d} :=\Phi^{-1}(F_{\nu_l}(\Hstat_{\theta_{l,d}}(T_{l,d}))),
\end{equation} 
where $1 \leq l \leq L$, $d \in \Dstat$. Observe that, by construction, if the marginal models are accurately specified and fitted, then $F_{\nu_l}(\Hstat_{\theta_{l,d}}(T_{l,d}))$, as a cdf, is uniformly distributed, and thus the $z_{l,d}$ variables from \eqref{eq:definez} would be standard normal distributions. This is the assumption we make for the rest of this section. 

With the marginals specified, the spatial dependence between the temperatures will be controlled by a Gaussian copula.
Specifically, we need to model the dependence of $z_{l,d}$ on the location $l$ and the day $d$, for $1 \leq l \leq L$, $d \in \Dstat$. We now assume that any fine (within each month) temporal dependence was accounted for with our seasonal modeling and that $z_{l,d_1}=z_{l,d_2}$ if $d_1$ and $d_2$ are daily indices belonging to the same month. Therefore, we combine all transformed monthly data at the same site and index it by month, resulting in about $30 \text{ (days per month)}\times 365 \text{ (days per year)} \times 29 \text{ (years)}$ realizations of $z_{l,m}$ for each $1 \leq m \leq 12$ month and $1 \leq l \leq L$ location. 

For each month $m$, the transformed temperature vector is modeled as a mean-zero Gaussian random vector with exponential correlation function $\text{Cor}(z_{m,l},z_{m,l'}) = e^{- \| \mathbf{s}_l - \mathbf{s}_{l'} \|/\gamma_m}$; the dependence between months is ignored. 
We thus estimate a separate length-scale parameter $\hat \gamma_m$ in the exponential correlation for each month of the year; once that is carried out, the model is complete. 
Since the marginal and copula steps described above estimate parameters with maximum likelihood separately, some statistical efficiency is lost. However, this is common practice in copula literature since the joint modeling, and thus, joint fitting would be exceedingly complex \citep{joe2005}.

\paragraph{Sampling procedure for normal temperature scenarios}
To perform an unconditional simulation for one year at each site $l$ on a given day $d$ at a given hour $h$, we proceed in reverse. That is, each day $d$ is mapped into a month $m$, and a scaled random variable $z_{l,d}$ is sampled from a Gaussian process with exponential correlation and parameter $\hat \gamma_m$.
A realization of the temperature $T_{l,d}$ is obtained by computing $ \left(F{\nu_l}(\Hstat_{\theta_{l,d}})\right)^{-1} \circ \Phi (z_{l,d})$. Here $\Phi$ is the cdf of the standard normal, $^{-1}$ denotes the functional inverse of the cdf, and $\circ$ is the functional composition.

\paragraph{Sampling procedure for extreme temperature scenarios}

\label{ref:extremesample}
We observe that the unconditional simulation from our statistical model amounts to simulating a Gaussian random vector $z_{l,d}$ and then applying a component-wise monotone transformation $ \left(F{\nu_l}(\Hstat_{\theta_{l,d}})\right)^{-1} \circ \Phi (z_{l,d})$ which transforms each component of the vector to follow the estimated marginal distribution. In particular, this implies that if we have two realizations of $z_{l,d}$, $z^1_{l,d}$ and $z^2_{l,d}$, at fixed day $d$ and hour $h$, and one of them is dominated by the other, that is 
$z^1_{l,d} \geq z^2_{l,d}$, $\forall 1 \leq l \leq L$, then it must be true that the same relationship will hold between the mapped temperatures $T^1_{l,d} \geq T^2_{l,d}$, $\forall 1 \leq l \leq L$ and thus the corresponding spatial averages $T^1_{a,d}:=\frac{1}{L} \sum_{l=1}^L T^1_{l,d}$, and $T^2_{a,d}:=\frac{1}{L} \sum_{l=1}^L T^2_{l,d}$. Therefore, $z_{l,d}$ with large spatial averages $z_{a,d}:=\frac{1}{L} \sum_{l=1}^L z_{l,d}$ are more likely to map to temperatures $T_{l,d}$ with large spatial averages $T_{a,d}$. 

We can then follow \citet{owen2019} to exploit the properties of this structure in our simulation of extreme temperature scenarios, which we define to be those where the spatial average is above its 0.99 quantile or below its 0.01 quantile.
The value for $T_{a,d}^{h,0.99}$, the high extreme quantile of the spatial average at day $d$ and hour $h$, is approximated with the corresponding quantile of the spatially averaged temperature from 100,000 unconditional Monte Carlo simulations from the temperature model at that $d$ and $h$ (and similarly for $T_{a,d}^{h,0.01}$). 
The next step is to produce temperature fields whose spatial average is more extreme than $T_{a,d}^{h,\pstat}$, for $\pstat=0.99$ and $\pstat=0.01$.
Here, brute-force accept-reject Monte Carlo simulation would be very computationally expensive due to the many rejections. 
Instead, since the distribution of a Gaussian random vector conditional on its average is available in closed form, we can easily produce samples of the scaled process $z_{l,d}$ at fixed day $d$ and hour $h$ conditional on the spatial average $z_{a,d}=\frac{1}{L} \sum_{l=1}^L z_{l,d}$ exceeding a prescribed threshold; for example the quantile of the standard normal $z_\alpha$ for $\alpha \leq 0.99$ but not too far smaller than it. We typically choose $\alpha=0.975$ for achieving $\pstat=0.99$ for temperature; that is, we sample conditional on $z_{a,d} \geq 1.96 = z_{0.975}$ and if the sample mapped into temperature $T_{l,d}$ has a spatial average less than $T_{a,d}^{h,0.99}$, we reject it. While this may introduce a bit of bias since it is conceivable that there may be samples with $z_{a,d} \leq 0.975$ where  $T_{a,d} \geq T_{a,d}^{h,0.99}$ the probability is insignificant for large $L$. 
This idea is implemented in a scheme to efficiently generate 100,000 scenarios when the spatial average temperature is extremely warm (two representative dates from mid-July and August) and extremely cool (two representative dates from mid-January and February).
These collections of extreme temperature simulations, as well as unconditional normal simulations, are used as input to the temperature dependence modeling scheme described in Section \ref{sec:temperature-dependence}.

\paragraph{Limitations of the temperature model}

We make several simplifying assumptions in our procedure, as modeling the entire spatio-temporal distribution of temperature while prioritizing extreme values is a  challenging task.
A more desirable spatial extreme model would enforce spatial continuity of the marginal BATs parameters, likely using a Bayesian hierarchical framework.
In our case, we manually corrected the shape parameters of the marginal BATs distribution since the estimated values were noticeably different at a few nearby counties.
Such differences in tail behaviors at nearby locations produced unrealistic temperature fields (e.g.,\ samples significantly outside the temperature record). 
Incorporating more sophisticated spatial dependence, as well as addressing the assumption of stochastic independence over days, could reduce these undesirable artifacts. 
The Gaussian copula \eqref{eq:gaussiancopula} can be restrictive for modeling joint extremes since it exhibits zero tail dependence \citep{MR1999654}; a student-$t$ copula may be a more realistic assumption. 
Despite these shortcomings, our spatial model serves as an efficient generator for temperatures in both normal and extreme scenarios.

\section{Case study: Midwest power grid planning}\label{sec:casestudy}
We conduct a capacity expansion case study focusing on the midwestern United States.  The operational subproblems are designed to capture key temperature dependencies in the power grid to evaluate the importance of such dependencies. Consistent with the model described in Section \ref{sec:modeling-approach}, we make a ``copper plate'' assumption (i.e., transmission constraints are not considered in any way) and do not model generator operation requirements such as ramping rates and commitment decisions. Due to these simplifying assumptions, the results of the case study should not be interpreted as making any actual recommendations for where to install generation capacity in any actual grid.
Instead, the purposes of this case study are to understand the relative merits of the proposed models for supporting decision-making in the presence of rare, high-impact events, to investigate the value of the proposed spatial temperature distribution model, and to investigate the importance of modeling temperature dependence of generation capacity when making capacity expansion decisions.


\subsection{Grid modeling assumptions}
Consistent with the decision variables introduced in Section \ref{sec:gen-decisions}, the generators available to build are of two varieties: (1)  combined cycle, combustion turbine, coal, nuclear, solar, and wind generators, which have binary variables and have an associated capital cost ($\gensbin$); and (2) potential new solar and wind for each county, which have nonnegative, bounded continuous variables and a per MW capital cost ($\genscont$). 

We use the fleet of existing generators as the set of candidate generators $\gensbin$ -- i.e., for each generator in the real fleet, we treat it as a ``candidate'' generator to be built or not. Information about the existing generators, including their efficiencies and locations, are from the EPA Needs database. The term ``generator'' here refers to a ``unit'' in the EPA data, so each ``unit'' at a larger generation facility is assigned its own decision variable. The assigned capital costs for these generators are based on EIA.gov estimates for the per megawatt cost of new construction of these plants. In the context of wind and solar power, we use the term ``generator'' to refer to a wind or solar farm. Wind and solar generator capital costs and capacities are NREL estimates on a county-by-county basis for new construction \citep{nrelreeds}. The wind and solar hourly capacity factors are also NREL estimates on a county-by-county basis \citep{nrelnaris}. All capital costs are amortized with a time horizon of 30 years to be on the same cost basis as the operating costs described below. The data described in this paragraph were initially collected and processed as part of the WEREWOLF project described in \citet{werewolf}.

Existing generators of other types (biomass, non-fossil waste, landfill gas, oil, and gas steam plants, hydropower) are included in the model and do not have first-stage (i.e. $\sone$) decision variables or cost ($\costcap$) associated with them, but are assumed to be available for dispatching in the second stage. The total capacity from these generators is less than 10\% of the total existing capacity. Additionally, to decrease the computational burden in running the experiments, the following decisions were fixed: half of the combustion turbines and half of the combined cycle units were chosen (i.e., their corresponding $\sone_i$ variables are fixed to $1$), and half of the coal units were excluded (i.e., their corresponding $\sone_i$ variables are fixed to $0$). Those choices were made based on preliminary experiments that indicated that in most optimal solutions the majority of the combined cycle units are selected and most of the coal units are not selected.

In addition to capital costs, each generator $i$ with a fuel source (except nuclear) has a per MW operating cost $\costgen_i$ based on their generator-specific efficiency and recent fuel costs (we used \$2.56 per MBtu for natural gas and \$2.02 per MBtu for coal). Nuclear is treated as free to operate if it is built. Hydro, solar, and wind also have no operating costs.

The cost of load shedding is $\costshed$ per MW. In the bi-objective models, this parameter is fixed to be slightly higher than the highest per unit operational generation cost among all candidate generators, as, in this model, the intent is just to have the cost high enough so that in the operation phase load is shed only if the demand exceeds total available capacity. As described in Section \ref{sec:basesolve}, in the base model, the parameter  $\costshed$ is varied to give solutions having different levels of risk of load shed.

Wind and solar base capacity factors are sampled according to the season and time of day of the sampled hour $H$ from a year's worth of hourly, county-specific capacity factor estimates from NREL described above.

\subsection{Scenario construction}
Here, we describe how samples of the parameters $\time,\temp,\outages,\demand$ are constructed, with further details found in the sections below. A ``normal'' scenario is one that uses an unconditional temperature distribution, while an ``extreme'' scenario is one that uses a conditionally extreme temperature distribution (as described in Section \ref{sec:temp-model}). 

For modeling ``normal'' scenarios, we use a set of hours that are chosen at random, uniformly among all even hours of the year and stratified among the four seasons. This set of hours is then fixed and declared as the set of hours $\mathcal{H}$ that are then used in the construction of all samples of normal scenarios. Each $H \in \mathcal{H}$  corresponds to a season and an hour $h$ of the day. Thus, for each $H \in \mathcal{H}$ some number of samples of the other uncertain parameters is taken to construct full scenarios. Specifically, the hour $H$ is used to sample a temperature scenario $\temp$ from the distribution defined for that particular hour (as described in Section \ref{sec:temp-model}). The hour of the year also defines the base wind and solar capacity factors (Section \ref{sec:temperature-dependence}). Finally,  the hour of the year and the simulated temperatures are used to generate an observation of the total electric load $\demand$ and the generator outages and de-ratings $\outages$ (also Section \ref{sec:temperature-dependence}). 

For ``extreme'' scenarios, the four seasons are replaced by four months (January, February, July, August) to focus on the coldest and hottest times of the year, and the remainder of the process is the same as for ``normal'' scenarios but with temperature samples being taken from the conditional distribution as described in Section \ref{sec:temp-model}.

\subsubsection{Temperature dependence of electricity generation}\label{sec:temperature-dependence}
The model of temperature dependence of generator outputs varies according to the generator technology type. The temperature dependencies described below are captured by the vector $\Outages$ of random variables, whose components $\Outages_i$ take on values greater than or equal to 0 (and usually less than or equal to 1) for each generator $i$.

Conditional on an observed outcome of the temperature vector, coal, gas, hydro, and nuclear plants may, randomly and independently of each other, experience a forced outage for a scenario (i.e. have capacity 0 for that hour). The probability of an outage is a function of the temperature in the plant's county and follows the technology-specific availability curves found in \citet{murphy2019}\footnote{The curves are found in Figure 6 of \citet{murphy2019} while the corresponding values are the ``Temperature-dependent forced outage rates'' in Table 5 of \citet{murphy2020}. In this paper, these values are linearly interpolated between the $5^\circ\celsius$  increments and linearly extrapolated beyond the temperature range using the slope between the hottest or coldest two values.}. That is, for $i \in \gensbin$, $\Outages_i $ is a Bernoulli random variable with probability equal to one minus the estimated (temperature-dependent) forced outage rate for generator $i$'s type for the temperature at generator $i$'s location.

In cold temperatures, wind farms may experience de-rating. Given a capacity $\sone_i$ and capacity factor $\cap_i(H)$, a wind farm $i$ has a random de-rating amount $\Outages_i \in [0,1]$ such that the output of the wind farm is $x_i p_i(H) \Outages_i$. Based on the temperature $\Temp$, the random variable $\Outages_i$ for a wind farm $i$ in location $l_i$ is defined as follows:
    \begin{equation}\label{eq:temp-depend-wind}
        \Outages_i \sim \begin{cases} \textbf{ 1}, & \Temp_{l_i} \geq -20^\circ \fahrenheit;\\ N(q, q(1-q)k), & -30^\circ \fahrenheit < \Temp_{l_i} < -20^\circ \fahrenheit;\\ \textbf{0}, & \Temp_{l_i} \leq -30^\circ \fahrenheit, \end{cases} 
    \end{equation}
where $N(q,q(1-q)k)$ indicates a normal distribution with mean $q$ and standard deviation $q(1-q)k$, $k=0.1$ and $q = \frac{0.1(T_{l_i}+20)}{ 31 + T_{l_i}}$. Here, \textbf{1} is defined as a random variable that always takes the value of 1, and similarly for \textbf{0}.
    
Photovoltaic solar panels decrease in efficiency as the panel temperature increases and increase in efficiency as panel temperature decreases \citep{landis1994}. For a solar farm $i$ at hour $H$ with a capacity factor $p_i(H)$, this relationship is modeled as
    \begin{equation}\label{eq:temp-depend-solar}
    \Outages_i = 1-(\Temp_{l_i} - 25^\circ\celsius + C_1\cdot p_i(H))\cdot \eta, 
    \end{equation}
for constants $C_1 = 25/144$ and $\eta=0.005$ and is derived as follows. The panel temperature is modeled using the following parameter assumptions: $T_\text{NOCT} = 45^\circ C$ is the normal operating conditions temperature of the panel, $e = 18\%$ is the efficiency of the panel, and $\eta = 0.005$ is the temperature coefficient. Given the capacity factor, total irradiance onto the panel is approximated by $p_i(H)/e$, so the temperature of the panel is, according to the Ross equation \citep{olukan2014},
    \begin{equation}\label{eq:ross}
    \Temp_{l_i,\text{cell}} = \Temp_{l_i} + \frac{T_\text{NOCT} - 20^\circ \celsius}{800 \text{W/m}^2}\cdot \frac{p_i(H)}{e}. 
    \end{equation}
Applying the temperature coefficient $\eta$, the temperature-dependent capacity factor is
    \begin{equation}\label{eq:solar}
    \Outages_i = 1 - (\Temp_{l_i,\text{cell}} - 25^\circ \celsius)\cdot \eta.
    \end{equation}
Substituting to match \eqref{eq:temp-depend-solar}, we find $C_1 = \frac{T_\text{NOCT} - 20^\circ \celsius}{e \cdot 800 \text{W/m}^2} = 25/144.$ Note that in this case, it is possible for $\Outages$ to take on values greater than 1 (e.g., in cold temperatures with high capacity factors). Additionally, this approach does not model impacts like snow cover that are not directly tied to temperature.

\subsubsection{Temperature dependence of electricity demand}
Electricity demand is also temperature-dependent. Demand is modeled using linear regression with temperature, the hour of the day, and season of the year (Winter: Dec, Jan, Feb; Spring: Mar, Apr, May; Summer: Jun, Jul, Aug; Fall: Sep, Oct, Nov) as the input data.
This linear regression is trained using the same temperature data as used in Section \ref{sec:temp-model} (NLDAS-2) and load data from the EPA Needs database.
Each county is modeled independently. Given an hour of the day $h$ and a corresponding temperature scenario $\temp$, the demand $\demand$ is the aggregate demand across all counties in the region.
To obtain the load per county, we build a linear regression model for each county to predict load based on temperature data in a manner described in \citet{ieee2019}. The following features are used to predict the load at hour $h$:
    the temperature $T_h$,
    the temperature squared $T_h^2$,
    an indicator $T_h^* = \indicator\bracket{T_h \geq 65^\circ \fahrenheit}$,
    the season $M_h$ of the year,
    and the day $W_h$ of week, where the options are weekend or weekday.
The parameters $T_h^*$, $M_h$, and $W_h$ are 1-hot encoded to generate the following features in the regression model: $M_h$, \ $T_h \cdot T_h^* \cdot M_h$, \ $T_h^2 \cdot T_h^* \cdot M_h$, \ $T_h \cdot T_h^* \cdot h$, \ $T_h^2 \cdot T_h^* \cdot h$, \ and $W_h \cdot h$.

Given the trained regression models (one per county), a total demand scenario is obtained by first determining the hour $H$ and the observation of the temperatures at all the counties. A demand scenario at each county is obtained by passing the hour of day $h$ and season of year (both determined by $H$) and the generated county temperature into the regression model for that county (for simplicity, all scenarios are considered weekdays in this step). The total demand in that sample is then obtained by summing the demands over all counties.

\section{Analyses}\label{sec:analyses}
We use the case study setup described in the previous section to perform experiments to assess the value of our approach. Section \ref{sec:experimental-details} describes the experimental design, and the remaining subsections describe specific analyses. These are the models considered in those analyses, which are summarized in Table \ref{tab:models}.
\begin{itemize}
    \item Base: The base model \eqref{eq:base-model-saa}.
    \item BO-$\cvar{}$: The bi-objective model \eqref{eq:biobj-unc-cvar}, with $\cvar{0.0001}$ as the risk measure.
    \item BO-Cond-$\cvar{}$: The bi-objective model \eqref{eq:biobj-cond-cvar}, with the conditional temperature samples described in Section \ref{sec:biobj-conditional} and $\cvar{0.1}$ as the risk measure.
    \item BO-Cond-$\cvar{}$ (Ind): The same as BO-Cond-$\cvar{}$, except the temperature samples are \textit{not} spatially dependent. That is, the correlations in temperature across space are ignored in both the conditional and unconditional samples of $\Temp$.
    \item BO-Cond-\lolp: The same as BO-Cond-$\cvar{}$, except the risk measure is \lolp\ instead of $\cvar{}$ (i.e., model \eqref{eq:biobj-cond-lolp})
    \item BO-Cond-$\cvar{\cvarthres}$: The same as BO-Cond-$\cvar{}$, except the threshold for $\cvar{}$ is $\cvarthres$.
    \item BO-Cond-$\cvar{}$ (NT): The same as BO-Cond-$\cvar{}$, except the generator temperature dependence modeling is removed. Instead, conventional generators (those in $\gensbin$) are modeled with an average forced outage rate (not one that changes as a function of $\temp$), and renewable generators (those in $\genscont$) have $\outages = 1$.
\end{itemize}
The $\cvar{}$ threshold used in the conditional model of $\cvarthres = 0.1$ is chosen by convenience: it is large enough to be less sensitive to sampling error while small enough to be sensitive to extremes. Section \ref{sec:riskmeasurecomp} demonstrates that the exact value of $\cvarthres$ is not that important for obtaining good solutions with this model. The temperature samples are conditional on the 0.01 and 0.99 quantiles (referred to as 1\% extreme).

\begin{table}[h]
    \centering
    \begin{tabular}{l|c l c c c}
Name & Model & \thead{Risk \\measure} & \thead{Condition \\ on $\objb$ \\ scenarios} & \thead{Spatially \\ dependent\\ $\temp$?} & \thead{Temperature \\ dependent \\ supply?} \\\hline
Base & \eqref{eq:base-model-saa} & None & N/A & Yes & Yes \\ 
BO-$\cvar{}$ & \eqref{eq:biobj-unc-cvar} & $\cvar{0.0001}$ & None & Yes & Yes \\ 
BO-Cond-$\cvar{}$ & \eqref{eq:biobj-cond-cvar} & $\cvar{0.1}$ & 1\% extreme & Yes & Yes \\ 
BO-Cond-$\cvar{}$ (Ind) & \eqref{eq:biobj-cond-cvar} & $\cvar{0.1}$ & 1\% extreme & \textbf{No} & Yes \\ 
BO-Cond-\lolp & \eqref{eq:biobj-cond-lolp} & $\lolp$ & 1\% extreme & Yes & Yes \\ 
BO-Cond-$\cvar{\cvarthres}$ & \eqref{eq:biobj-cond-cvar} & $\cvar{\cvarthres}$ & 1\% extreme & Yes & Yes \\ 
BO-Cond-$\cvar{}$ (NT) & \eqref{eq:biobj-cond-cvar} & $\cvar{0.1}$ & 1\% extreme & Yes & \textbf{No} 
    \end{tabular}
    \caption{Models considered in the analyses.
    }
    \label{tab:models}
\end{table}

\subsection{Experimental details}\label{sec:experimental-details}
All models are implemented in Python using the extensive forms described in the Appendix \ref{app:forms}. These formulations are solved with the MIP solver Gurobi using the HTCondor system from the Center for High Throughout Computing at UW-Madison \citep{chtc}. This system is a collection of thousands of shared, multi-core computers to which users can submit computational jobs.

\subsubsection{Solution generation}\label{sec:training}
Some small variations are mentioned in later subsections, but the following applies to all models. To apply SAA, the models are solved over 15 samples each with 1056 scenarios. To find solutions with different trade-offs between cost and risk, the models are solved for multiple values (roughly 10-15) of the parameter that controls risk preference. For the base model, the parameter is $\costshed$, while for the bi-objective models the parameter is the bound on the second objective in the $\epsilon$-constraint method, described in Section \ref{sec:biobj-unconditional}. The specific values of these parameters used are found in Appendix \ref{appendix:parameter-values}. This means results from one experiment for one model solved with 15 samples and 10 parameter values may yield up to 150 different solutions, although some samples and parameter values may yield the same solution, so  there may be less than 150 distinct solutions.

When generating solutions using the base model each SAA uses a sample of 1056 scenarios. For the bi-objective models each SAA uses 528 scenarios for the first objective and a different 528 scenarios for the second objective. (Thus, the total number of scenarios used in the different methods is held constant at 1056.) The models were terminated when they reached the optimality gap tolerance of 0.01\% or a time limit of 6 hours. Each instance was given access to up to 160 GB of memory. Most instances were solved within the optimality gap tolerance in much less than these limits: a typical BO-Cond-$\cvar{}$ instance solved in about 10 minutes and required about 25 GB of memory. While the heterogenous nature of the machines used in HTCondor makes any comparison very approximate, and hence this is not a focus of our study, some significant differences in solution time were observed: using \lolp\ as the risk measure or solving the base model occasionally required resources closer to the memory limit or hit the time limit.

\subsubsection{Solution evaluation}\label{sec:testing}
All of the solutions are tested using a set of evaluation scenarios that is different from the scenarios used in the SAA problems used to generate the solutions. The set of evaluation scenarios is the same for all methods. Testing a solution involves using Monte Carlo estimation to evaluate its performance on both objectives of our bi-objective model over larger samples. One large sample (25,008 scenarios) is used to estimate the second objective (extreme risk). These are unconditional scenarios generated the same way as the other unconditional scenarios and are intended as an approximation of the ``true'' temperature distribution. For each scenario $k$, the extreme subproblem $\objb$ is evaluated by \eqref{eq:qextrememax}, yielding $\shed_k$ of load shed. The \lolp\ of the solution is then estimated as $\frac{1}{25008}\sum_{k=1}^{25008} \indicator\bracket{\shed_k > 0}$. The $\cvar{\cvarthres}$ of the solution is estimated by computing the mean of the $\lceil25008 \cdot \cvarthres\rceil$ largest $\shed_k$ values. For our choice of $\cvarthres = 0.0001$, this means that we take the mean of the three highest load shed values.

We use a smaller subset (5040 scenarios) of the large evaluation sample to estimate the first objective (expected cost). For each scenario $k$, the normal subproblem $\obja$ is evaluated by solving the linear program \eqref{eq:qnormal} with a fixed value of $\costshed$, which yields a cost $t_k$. The same value of $\costshed$ is used here as is used when solving the bi-objective models (i.e., slightly larger than the largest marginal cost of generation from any generator). The total cost of a solution is estimated by $\frac{1}{5040} \sum_{k=1}^{5040} t_k$.

\subsubsection{Trade-off curve}\label{sec:tradeoffcurve}
As argued in Section \ref{sec:biobj-model}, the decision maker must choose a solution with a suitable trade-off between risk and cost. Therefore, our goal is to generate solutions along the whole efficient frontier, or trade-off curve. To assess methods with this goal in mind, we plot the trade-off curve of the solutions obtained by each method, which corresponds to a scatter plot in which the $x$-axis is an estimate of the measure of risk (either \lolp\ or $\cvar{}$), and the $y$-axis is an estimate of the expected total cost. In some figures, we exclude solutions dominated by other solutions found by the same method. For instance, in an experiment, if the base model finds a solution that has higher estimated risk and expected cost than another solution found by the base model, it is a dominated solution and would not be shown in the figure.

We choose a $\cvar{0.0001}$ value of 15 GWh as a cutoff for consideration and, in most figures, do not display solutions having risk levels higher than this. This value is approximately the point at which solutions have significant load shed in the test scenarios for the first objective, indicating they would be unacceptable to any decision-maker. 

\subsection{Base versus bi-objective (conditional \& unconditional)}\label{sec:base-vs-biobj}
\begin{figure}[h]
    \centering
\includegraphics{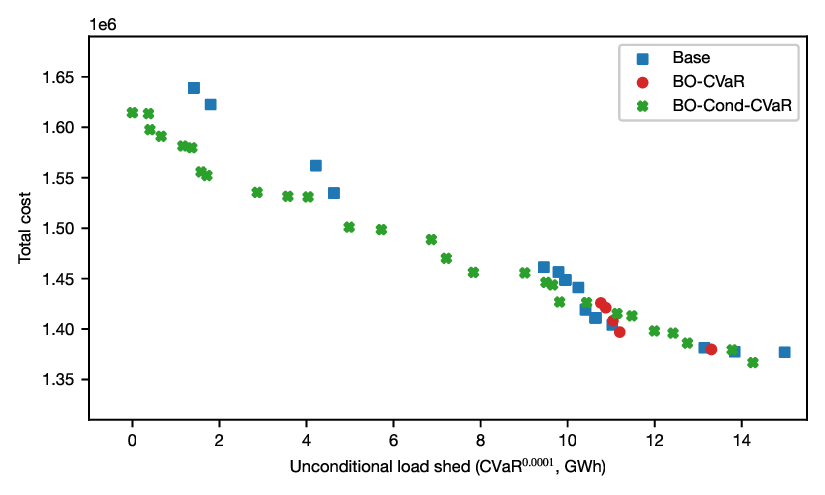}
\caption{Trade-off curve showing objective values for solutions from Base, BO-$\cvar{}$ and BO-Cond-$\cvar{}$. Dominated solutions and solutions with $x$-axis value greater than 15GWh are not shown.}
\label{fig:trad-off_base-biobj-cond}
\end{figure}

\begin{figure}[h]
    \centering
\begin{subfigure}{0.45 \textwidth}
\includegraphics{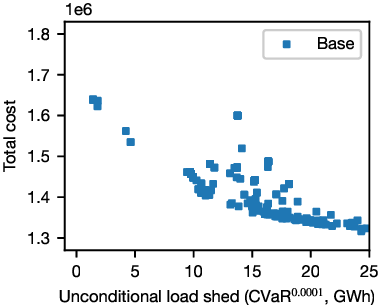}
\caption{Base}
\label{fig:small-trade-off_base}
\end{subfigure}
\begin{subfigure}{0.45\textwidth}
\includegraphics{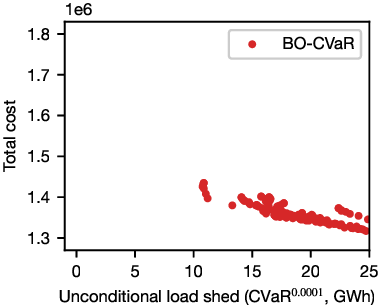}
\caption{BO-$\cvar{}$}
\label{fig:small-trade-off_biobj}
\end{subfigure}
\begin{subfigure}{0.45\textwidth}
\includegraphics{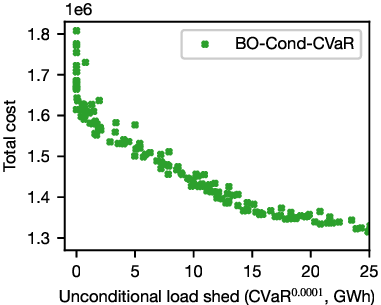}
\caption{BO-Cond-$\cvar{}$}
\label{fig:small-trade-off_cond}
\end{subfigure}
\begin{subfigure}{0.45\textwidth}
\includegraphics{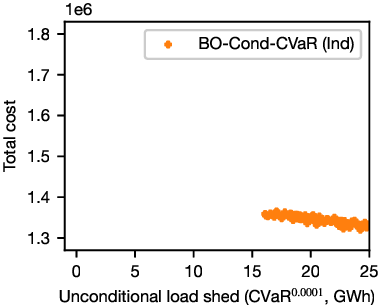}
\caption{BO-Cond-$\cvar{}$ (Ind)}
\label{fig:small-trade-off_ind}
\end{subfigure}

\caption{Trade-off curves for four models with dominated solutions included.}
\label{fig:small-trade-off}
\end{figure}

We first investigate the abilities of Base, BO-$\cvar{}$ and BO-Cond-$\cvar{}$ at finding solutions along the efficient frontier. After solving each model over each sample and each parameter value, all solutions are evaluated as described in Section \ref{sec:testing}. The trade-off curves showing these results are shown in Figure \ref{fig:trad-off_base-biobj-cond}, where the non-dominated solutions for each method are plotted, and Figures \ref{fig:small-trade-off} (a)\---(c), where dominated solutions are included. BO-Cond-$\cvar{}$ produces the most filled-in trade-off curve, which includes solutions with no load shed. Base finds solutions at different points along the curve, including occasional solutions with low risk, but has large gaps between points along the curve, and its lower risk solutions are dominated by BO-Cond-$\cvar{}$. BO-$\cvar{}$ finds efficient solutions at high risk but does not find lower-risk solutions.

These results show that while all methods can find non-dominated solutions, BO-Cond-$\cvar{}$ does a better job of finding solutions along the full range of the trade-off curve, including solutions with very low risk. Neither Base nor BO-$\cvar{}$ trace out as complete of a trade-off curve compared with BO-Cond-$\cvar{}$, which would hinder a decision maker's ability to choose a preferred trade-off. Figures \ref{fig:small-trade-off} (a)\---(c) demonstrate that BO-Cond-$\cvar{}$ consistently finds solutions close to the efficient frontier across the multiple samples. In contrast, the BO-$\cvar{}$ and more significantly Base yield more solutions that are farther from the efficient frontier. This suggests that BO-Cond-$\cvar{}$ is less sensitive to sampling error than the alternatives.

\subsection{Why does conditional approximation work?}\label{sec:conditionalapprox} 
Proposition \ref{proposition:pareto-monotone} in  Section \ref{sec:biobj-conditional} indicates that if $\cvar{0.1}\bracket{\objb\paren{\sone,\RandE^{\Ntwo^E}}}$ (the second objective of BO-Cond-$\cvar{}$) is approximately equal to some monotone transformation of \allowbreak $\cvar{0.0001}\bracket{\objb\paren{\sone,\Rand^{\Ntwo}}}$ (the second objective of BO-$\cvar{}$), then we can expect the set of Pareto solutions obtained using BO-Cond-$\cvar{}$ to be a good approximation of the set of Pareto solutions of BO-$\cvar{}$. 

\begin{figure}[h]
    \centering
\includegraphics{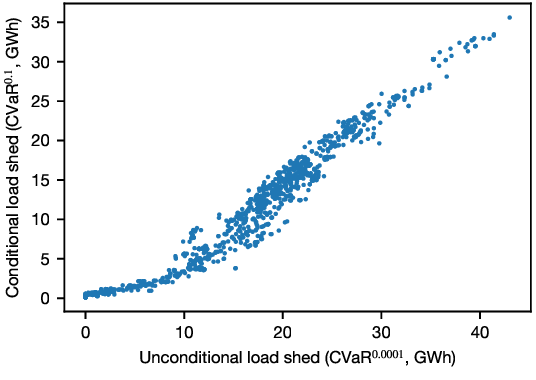}
\caption{Relationship between unconditional risk measure and conditional risk measure. The solutions with a total cost of less than $10^7$ have been plotted. The $x$-axis shows the unconditional risk measure value and the $y$-axis shows the conditional risk measure value. One additional solution, at $(89, 64)$, has been removed for presentation, though it is consistent with the increasing pattern.}
\label{fig:monotone-relationship}
\end{figure}

We thus conducted an experiment to test if $\cvar{0.1}\bracket{\objb\paren{\sone,\RandE^{\Ntwo^E}}}$is approximately equal to some monotone transformation of \allowbreak $\cvar{0.0001}\bracket{\objb\paren{\sone,\Rand^{\Ntwo}}}$ over a range of solutions $x$. Specifically, we evaluate the unconditional risk and conditional risk of all the first-stage solutions $x$ found in the experiment in Section \ref{sec:base-vs-biobj} as well as all solutions found by Gurobi during its solution process. The unconditional risk is evaluated via the standard evaluation procedure described in Section \ref{sec:testing}. The conditional risk is evaluated by using conditional scenarios instead of unconditional scenarios in the testing procedure for the second objective. These conditional scenarios (a total of 15,840 of them) include the samples used to train BO-Cond-$\cvar{}$. We exclude solutions that have extremely high costs on the first, cost objective (greater than $1e7$, while Pareto solutions are less than $1.7e6$) since these are clearly not relevant for approximating the set of Pareto solutions. We observe in Figure \ref{fig:monotone-relationship} that there is an approximately monotone relationship between the conditional and unconditional risk measures, indicating that the good performance of BO-Cond-$\cvar{}$ discussed in Section \ref{sec:base-vs-biobj} may be explained by Proposition \ref{proposition:pareto-monotone}. 

\subsection{Spatial distribution}\label{sec:comparespatialvsind}
We next investigate the value of using the spatially dependent temperature scenarios described in Section \ref{sec:stats-model}. We attempt to construct an efficient frontier using a temperature distribution that ignores the spatial dependence and compare the results to our previous results. For this comparison, we use the BO-Cond-$\cvar{}$ method for constructing the efficient frontier, since we found it works uniformly across a range of risk values. The method for approximating the efficient frontier with the spatially independent distribution is  identical to the BO-Cond-$\cvar{}$ method, with the only change being that we sample the scenarios from a temperature distribution that assumes spatially \textit{independent} temperature realizations. Once the solutions are obtained, we evaluate them using the scenarios drawn from the spatially \textit{dependent} model.  The trade-off curves for these models are shown in Figure  \ref{fig:small-trade-off_ind}. BO-Cond-$\cvar{}$ (Ind) finds efficient solutions at high risk consistently (with cost and risk indistinguishable from BO-Cond-$\cvar{}$) but does not find any low-risk solutions. This indicates that modeling spatial temperature dependence is essential for finding solutions with low risk.

\begin{figure}[h]
    \centering
\begin{subfigure}{0.45\textwidth}
\includegraphics{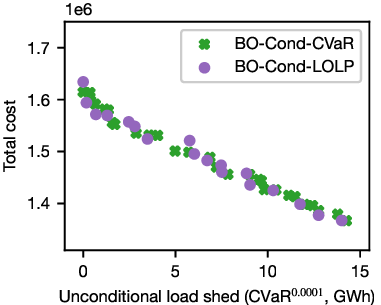}
\caption{Tested on CVaR}
\label{fig:risk-measures_cvar}
\end{subfigure}
\begin{subfigure}{0.45\textwidth}
\includegraphics{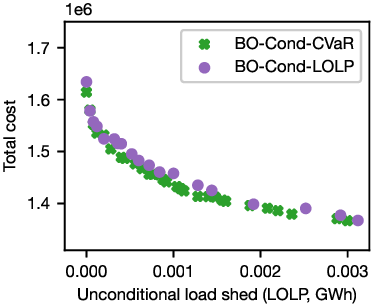}
\caption{Tested on LOLP}
\label{fig:risk-measures_lolp}
\end{subfigure}

\caption{CVaR vs LOLP. Both figures show the same solutions from each model and have cost on the $y$-axis. The left figure (a) evaluates $\cvar{0.0001}$ for the $x$-axis, and the right figure (b) evaluates \lolp\ for the $x$-axis.}
\label{fig:risk-measures}
\end{figure}

\subsection{Risk measure comparison}\label{sec:riskmeasurecomp}
We next compare the two risk measures (\lolp\ and $\cvar{}$) used to assess extreme risk. Our hypothesis is that the two measures will lead to similar sets of solutions on the trade-off curve. To test this, we find two sets of solutions, one using each metric, and evaluate them on the two  metrics to see how they compare.  The results, shown in Figure \ref{fig:risk-measures}, show very comparable performance between the two measures. When evaluated on $\cvar{0.0001}$, the solutions obtained using \lolp\ and $\cvar{}$ in the optimization models traced out very similar curves. Perhaps surprisingly, when evaluated on \lolp, the solutions obtained from the model optimizing $\cvar{}$\ slightly outperform those obtained from the model optimizing \lolp. The explanation for this slight disparity is that the \lolp\ optimization model is more difficult to solve, and hence occasionally was terminated due to time or memory limits, while this did not occur with the $\cvar{}$ model. Given that $\cvar{}$ is computationally easier in optimization than \lolp, these results suggest a preference for using $\cvar{}$ in the optimization step, even if \lolp\ is the  risk measure to be used in the evaluation step.

\begin{figure}[h]
    \centering
\includegraphics{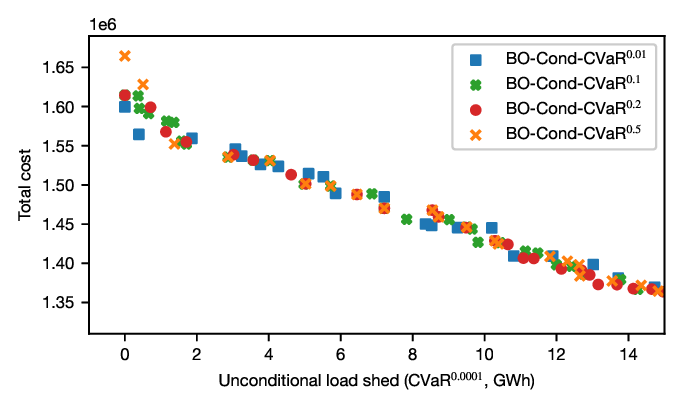}
\caption{Results from BO-Cond-$\cvar{\cvarthres}$ trained with different values of $\cvarthres$.}
\label{fig:cvar-comparison}
\end{figure}

We also investigate the sensitivity of $\cvar{\cvarthres}$ to the choice of $\cvarthres$ in BO-Cond-$\cvar{}$. As described in Section \ref{sec:biobj-conditional}, using conditionally extreme samples requires a different choice of $\cvarthres$, but the exact choice of $\cvarthres$ may be difficult. We thus explore the sensitivity of our results to this choice. In addition to $\cvarthres=0.1$ (the only value used so far), we solve BO-Cond-$\cvar{\cvarthres}$ with $\cvarthres \in \{0.01,0.2,0.5\}$. We then evaluate these solutions with the same unconditional testing scenarios described in Section \ref{sec:testing}. The non-dominated solutions are plotted in Figure \ref{fig:cvar-comparison} and indicate that the four choices of $\cvarthres$ produce similar trade-off curves. This suggests that these are all good choices of $\cvarthres$ and hence the precise choice $\cvarthres$ is not important.

\subsection{Generator temperature dependence}
In this section, we evaluate the impacts of modeling temperature dependence of generator availability. As noted in \citet{murphy2019}, a common practice when measuring a system's risk of load shedding is to assume generators' outages are not correlated with temperature; instead, it is assumed each generator has a ``forced outage rate,'' which is the probability that the generator will be unavailable. We design experiments to understand the impact of such an assumption compared to the temperature dependence of generator availability that has been used so far in our model.

\subsubsection{Base case study}\label{sec:temp-dependence-base}
We first use a similar approach as described in the previous section: we solve BO-Cond-$\cvar{}$ using different data. In this case, the new data will be identical to the ``correct'' data (the data used in Section \ref{sec:base-vs-biobj}) except for the following difference: the temperature dependence described in Section \ref{sec:temperature-dependence} is replaced by the following models of generator uncertainty $\Outages$:
\begin{itemize}
    \item Conventional generators are modeled with a temperature-independent forced outage rate. That is, for conventional $i \in \gens$, $\Outages_i$ is a Bernoulli random variable with probability equal to 1 minus generator $i$'s type's forced outage rate\footnote{The values used here are the ``Unconditional forced outage rates'' from Table 5 in \citep{murphy2020}.}.
    \item Solar and wind do not have any capacity adjustments due to temperature, so, for wind and solar generators $i \in \gens$, $\Outages_i \sim \mathbf{1}$.
\end{itemize}
Recall that the models of temperature dependence for solar generators increased their capacity factors in colder temperatures and decreased their capacity factors in higher temperatures. Therefore, removing this model leaves the mean capacity factor roughly the same, as is the effect for conventional generators. For wind generators, the prior temperature dependence model only affected output at temperatures below $-20^\circ\fahrenheit$. As a result, this model only affects performance in extreme temperatures, and the mean performance of wind is adjusted very little by this change.

\begin{figure}[H]
    \centering
\includegraphics{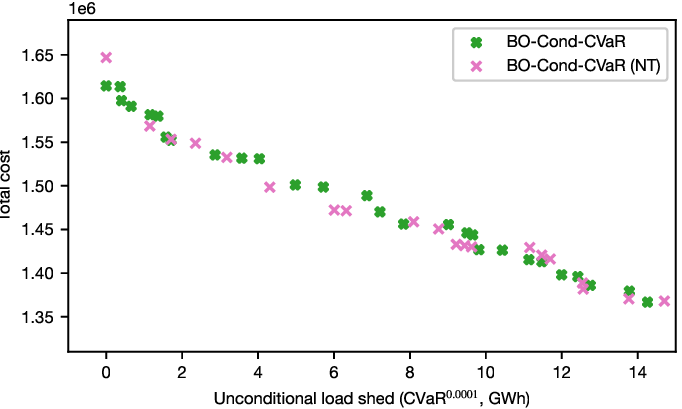}
\caption{Effects of modeling generator temperature dependence in base case study. BO-Cond-$\cvar{}$ corresponds to the results from Section \ref{sec:base-vs-biobj} and BO-Cond-$\cvar{}$ (NT) corresponds to solutions from the data with no temperature dependence of generator outages.}
\label{fig:notemp-comparison}
\end{figure}

The trade-off curves from these models can be seen in Figure \ref{fig:notemp-comparison}, which shows that both models find good solutions consistently at many different trade-offs between risk and cost. 
These results alone do not show a difference in performance between the models, suggesting that modeling generator temperature dependence is not necessary to finding efficient solutions at many risk levels in this case study.
Looking more closely though, the results demonstrate a degree of robustness of the bi-objective approach. Figure \ref{fig:notemp-by-U} in Appendix \ref{app:fig} displays all solutions found when solving both models using the $\epsilon$-constraint method, where the second objective (minimizing risk) is constrained to be at most $\bound$ while the first objective (average cost) is minimized.
For each choice of $\bound$, the solutions found by BO-Cond-$\cvar{}$ (NT) are, in general, higher risk and lower cost than the solutions found by BO-Cond-$\cvar{}$. This shift is compensated by solving the bi-objective model over many values of $\bound$.

Since our generator outage modeling considers the temperatures at each generator's location, the location of a generator has the potential to affect its desirability. For instance, between two identical generators, there is a preference for the one that is in a milder climate (and so experiences fewer outages). However, such a preference only exists when considering the generators' temperature-related outages. While the effects of not considering those outages were minimal in this experiment, we observed that the data of the case study limited the model's ability to choose generators based on their location: among ``reasonable'' solutions\footnote{Here, we mean those with $\cvar{0.0001}(\text{load shed}) \leq 15$ GWh when evaluated by the method in Section \ref{sec:testing}.}, many of the generators available in the model were always chosen because of the need to cover a large amount of load in every scenario. That left relatively few generators available to choose to reduce risk, and among those remaining, there existed large cost differences and little variety in location.  As a result, the full impact of generator temperature dependence was not able to be evaluated by this experiment.

\subsubsection{Modified case study}\label{sec:temp-dependence-modified}
In order to more effectively evaluate the potential impact of modeling generator temperature dependence, we modify the case study in the following ways:
\begin{enumerate}
    \item Additional new conventional generators are available in a quarter of the counties, with identical specifications across counties. The generator options are two sizes of gas combustion turbines, two sizes of gas combined cycle plants, and one size of nuclear plant. The full specifications are shown in Table \ref{tab:new-generators} in Appendix \ref{app:gens}.
    \item The fixed capital cost of existing conventional generators is reduced by 90\% (the new generators just mentioned take on the un-reduced cost per megawatt of installed capacity).
    \item Demand is increased by 25\%.
\end{enumerate}
These adjustments give the model the ability to make decisions based on generator location. Reducing the capital cost of existing generators reduces the significance of the variance among the less desirable generators and reflects a more realistic capacity expansion problem. Increasing demand requires that some new conventional generators be built instead of just relying on the existing fleet, and generators of identical specifications are available in a wide range of locations, exposing them to different temperature distributions.

\begin{figure}[H]
    \centering
\includegraphics{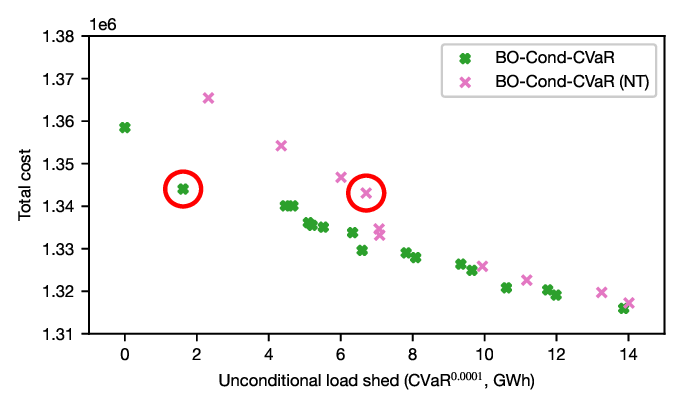}
\caption{Trade-off curves with nondominated solutions only of modified case study described in Section \ref{sec:temp-dependence-modified}. The model solved with generator temperature dependence (BO-Cond-$\cvar{}$) outperforms the model solved with no generator temperature dependence (BO-Cond-$\cvar{}$ (NT)) at lower-risk solutions. All solutions are evaluated on the same data that models generator temperature dependence.}
\label{fig:temp-modified}
\end{figure}

We repeat the same test on this modified model as we did in the previous section: We solve one model with the ``correct'' data (BO-Cond-$\cvar{}$) and one model with data that does not model generator temperature dependence (BO-Cond-$\cvar{}$ (NT)). We then evaluate them on an updated test set that models generator temperature dependence (the only changes from the original test set are to increase demand by 25\% and model generator outages for the additional generators in the same way as existing generators). The trade-off curves of this experiment are shown in Figure \ref{fig:temp-modified}. For higher risk values, the models perform similarly. However, BO-Cond-$\cvar{}$ (NT) has higher costs than BO-Cond-$\cvar{}$ for solutions with $\cvar{0.0001}$ values below 7 GWh.

\begin{figure}[h]
    \centering
\includegraphics{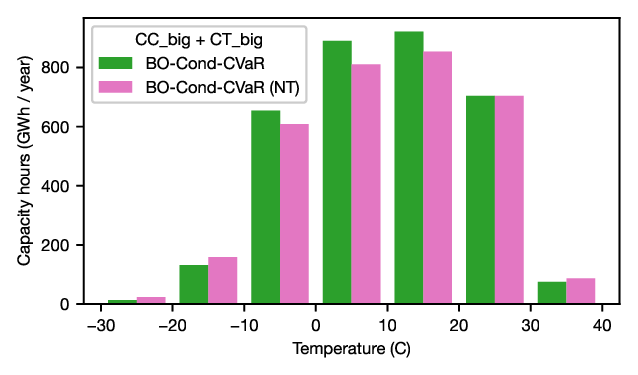}
\caption{Distribution of capacity-weighted new gas generation temperature exposure. }
\label{fig:histogram_gas_new-model}
\end{figure}

We next investigate how the solutions obtained by the model considering temperature dependence are able to achieve such lower risk at similar costs. To do so, we choose two solutions, one generated from each model, having similar expected cost; the selected solutions are circled in Figure \ref{fig:temp-modified}.  In Figure \ref{fig:histogram_gas_new-model} we display a histogram of the gas generation temperature exposure of the two selected solutions.
In this histogram, each new gas generator selected contributes its capacity to the corresponding temperature buckets based on the observed temperatures in the extreme objective test set. The total contribution of all generators is scaled to be yearly (i.e., multiplied by $8760 / |S|$, where $|S| = 25008$ is the size of the extreme objective test set).
This plot indicates that the solution selected by
BO-Cond-$\cvar{}$ includes generators that tend to experience less extreme temperatures. While not displayed in the figure, we also find that the solution from BO-Cond-$\cvar{}$ selects more combined cycle generators and fewer combustion turbines than the solution from BO-Cond-$\cvar{}$ (NT), which is consistent with the combustion turbine's more severe temperature dependence. 

\begin{figure}[h]
    \centering
    \includegraphics[width=5.5in]{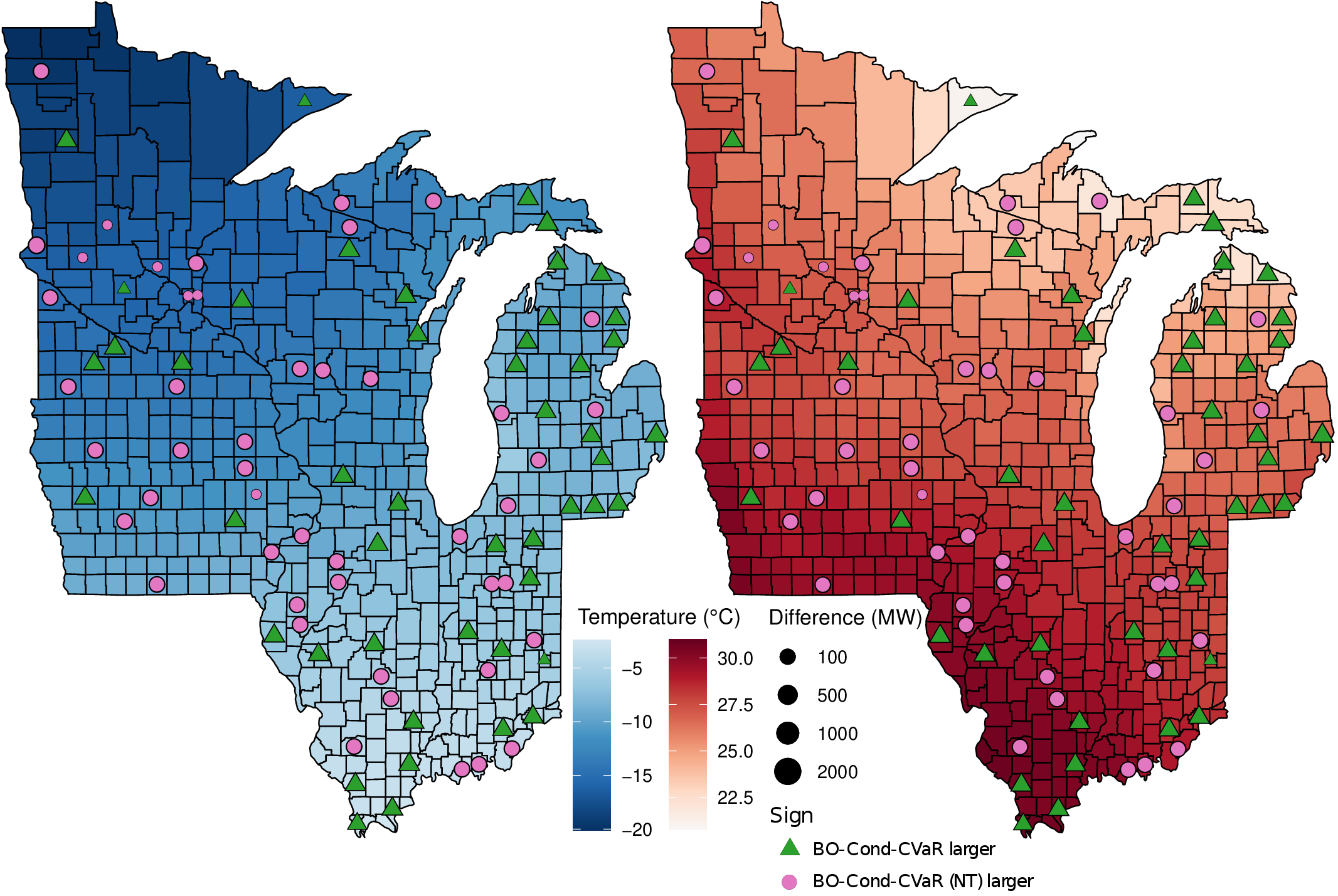}
    \caption{Geographic distribution of installed capacity in two selected solutions.}
    \label{fig:heatmap}
\end{figure}

We also look at the geographic distribution of generators in the two selected solutions. Figure \ref{fig:heatmap} shows the differences in new gas generation capacity along with the 0.05 and 0.95 temperature percentiles on a county-by-county basis.
This figure includes two maps of the region of interest showing each county. The blue shading in the left map shows each county's 0.05 percentile temperature (darker is colder). The red shading in the right map shows each county's 0.95 percentile temperature (darker is hotter). The circles and triangles show the county-by-county difference in new gas generation capacity between the two selected solutions. A green triangle indicates that the temperature dependent model (BO-Cond-$\cvar{}$) built more gas generation in that county, while a pink circle indicates that the non-temperature dependent model (BO-Cond-$\cvar{}$ (NT)) built more. The size of the triangle/circle indicates the magnitude of difference.
In general, the new gas generation chosen by BO-Cond-$\cvar{}$ is further east than what is chosen by BO-Cond-$\cvar{}$ (NT). This corresponds with less extreme cold temperatures in the southeast of the region and less extreme hot temperatures in the northeast of the region, suggesting BO-Cond-$\cvar{}$'s preference for generators that experience less extreme temperatures.
These observations demonstrate that generator temperature dependence modeling can make an impact the model's ability to find efficient, low risk solutions.

\section{Conclusion}
While our case study demonstrates the potential value of the proposed modeling approach, extensions of the model are needed before it can be useful for generator capacity expansion planning in practice.
Modeling these extensions, some of which are discussed in Section \ref{sec:model-extensions}, will increase the computational challenge of solving the resulting model and hence will likely require investigation of suitable solution techniques.

While this paper was focused on the problem of generation capacity expansion in the power grid, the challenge of designing a system that is efficient in ``normal'' operations but also protected against bad outcomes in rare, high-impact ``extreme'' conditions appears in many contexts. It would be interesting to explore whether our bi-objective modeling framework and conditional sampling approach could be useful in other such problems.  A key question in such exploration would be to 
determine criteria for defining ``extreme'' conditions, and to devise efficient methods for sampling conditional on being in such a condition.

\section*{Acknowledgements}
    This material was based upon work
supported by the U.S. Department of Energy, Office of Science,
Office of Advanced Scientific Computing Research (ASCR) under
Contract DE-AC02-06CH11347.
\clearpage
\bibliographystyle{apalike} 
\bibliography{bibl}

\appendix
\section{Appendix}
\subsection{Explicit extended formulations}
\label{app:forms}

\subsubsection*{Base model}
The extensive form of the base model \eqref{eq:base-model-saa}:
\begin{subequations}
    \begin{alignat*}{2}
        \min_{\sone \in X} \ & \costcap^\top \sone + \frac{1}{|N|}\sum_{k\in N} \paren{\sum_{i \in \gens} \costgen\paren{\stwo_{ik}} + \costshed \shed_k }\quad &&\\
        \st \ & \sum_{i \in \gens} \stwo_{ik} \geq \demand_k - \shed_k \quad && \forall k \in \N\\
        & \stwo_{ik} \leq \sone_{i} \outages_{ik} \cap_i\paren{\time_k} && \forall i \in \gens, \forall k\in \N\\
        & \stwo, \shed \geq 0.
    \end{alignat*}
\end{subequations}

\subsubsection*{Bi-objective model}
The extensive form of the bi-objective model with $\cvar{}$ as the risk measure using the $\epsilon$-constraint method, model \eqref{eq:biobj-unc-cvar-epsilon}:
\begin{subequations}\label{biobj:extensivecvar}
    \begin{alignat}{2}
        \min_{\sone \in X} \ & \costcap^\top \sone + \frac{1}{|\None|}\sum_{k\in \None} \paren{\sum_{i \in \gens} \costgen\paren{\stwo_{ik}} + \costshed \shed_k }\quad &&\\
        \st
        & \sum_{i \in \gens} \stwo_{ik} \geq \demand_k - \shed_k \quad && \forall k \in \None \nonumber\\
        & \stwo_{ik} \leq \sone_{i} \outages_{ik} \cap_i\paren{\time_k} && \forall i \in \gens, \forall k \in \None \nonumber\\
        & \stwo_k, \shed_k \geq 0  && \forall k \in \None \nonumber\\
        &\cvarvar + \frac{1}{1-\cvarthres} \sum_{k\in \Ntwo} \frac{\cvarvark_k}{|\Ntwo|} \leq \bound \quad && \label{eq:app:cvar1}\\
        &\cvarvark_k \geq \demand_k - \sum_{i \in \gens} \sone_i \outages_{ik} \cap_i \paren{\time_k} - \cvarvar, \quad &&\forall k\in \Ntwo \label{eq:app:cvar2}\\
        &\cvarvark_k \geq 0, \quad &&\forall k\in \Ntwo. \label{eq:app:cvar3}\\
        &\cvarvar \text{ unrestricted in sign} && \label{eq:app:cvar4}
    \end{alignat}
\end{subequations}
The constraints \eqref{eq:app:cvar1}-\eqref{eq:app:cvar4} are a formulation for $\cvar{\cvarthres}\bracket{\objb\paren{\sone,\Rand^{\Ntwo}}}$ due to  \citet{rockafellar2000}.

The $\epsilon$-constraint method applied to the bi-objective model with \lolp\ as the risk measure, model \eqref{eq:biobj-unc-lolp}:
\begin{subequations}\label{eq:app:lolp-epsilon}
\begin{align}
    \min_{\sone \in X} \ &  \costcap^\top \sone + \frac{1}{|\None|}\sum_{k\in \None} \obja \paren{\sone,\rand^k}\\
    \st & \frac{1}{|\Ntwo|} \sum_{k\in \Ntwo} \indicator \bracket{\objb\paren{\sone,\rand^k} > 0} \leq \bound. \label{eq:app:lolp-epsilon:constraint}
\end{align}
\end{subequations}
The extensive form of the above model \eqref{eq:app:lolp-epsilon}:
\begin{subequations}\label{eq:app:lolp-extensive}
    \begin{alignat}{2}
        \min_{\sone \in X} \ & \costcap^\top \sone + \frac{1}{|\None|}\sum_{k\in \None} \paren{\sum_{i \in \gens} \costgen\paren{\stwo_{ik}} + \costshed \shed_k }\quad &&\\
        \st \ 
        & \sum_{i \in \gens} \stwo_{ik} \geq \demand_k - \shed_k \quad && \forall k \in \None\\
        & \stwo_{ik} \leq \sone_{i} \outages_{ik} \cap_i\paren{\time_k} && \forall i \in \gens, \forall k \in \None\\
        & \stwo_k, \shed_k \geq 0  && \forall k \in \None\\
        &\sum_{k \in \Ntwo} \lolpvark_k \leq \left\lfloor |\Ntwo| \bound \right\rfloor &&  \label{eq:app:lolp-extensive:prob}
        \\
        & \bigM_k \lolpvark_k \geq \demand_k - \sum_{i \in \gens} \sone_i \outages_{ik} \cap_i \paren{\time_k} && \forall k \in \Ntwo \label{eq:app:lolp-extensive:bigM}
        \\
        & \lolpvark_k \in \{0,1\} && \forall k \in \Ntwo 
        \label{eq:app:lolp-extensive:chances}
    \end{alignat}
\end{subequations}
The chance constraint \eqref{eq:app:lolp-epsilon:constraint}, which limits the probability of load shed to be smaller than $\bound$, is formulated by constraints \eqref{eq:app:lolp-extensive:prob}, \eqref{eq:app:lolp-extensive:bigM}, and \eqref{eq:app:lolp-extensive:chances} in the extensive form: for $k \in \Ntwo$, the binary variable $\lolpvark_k$ takes on the value 1 if there is load shed in scenario $k$ (value 0 otherwise), and $\bigM_k$ is a value large enough so that constraint \eqref{eq:app:lolp-extensive:bigM} is not binding when $\lolpvark_k = 1$. 

We determine the values $\bigM_k$ by solving the following problem for all scenarios $k \in \Ntwo$:
\begin{subequations}\label{bigMproblem}
\begin{alignat}{3}
V_k = \min\ &\sum_{i \in \gens} \sone_i \outages^k_i \cap_i \paren{\time^k},\quad & \\
\st         & \sum_{i \in \gens} \sone_i \bar{\cap}_i \geq \demand_{1-U},   \quad & \\
            & \sone_i \in \{0,1\},                                          \quad & \forall i \in \gensbin, \\
            & 0 \leq \sone_i \leq \outages_i^k \cap_i\paren{\time^k} ,      \quad & \forall i \in \genscont,
\end{alignat}
\end{subequations}
where $\demand_{1-U}$ is the $1-U$ quantile of $\braces{ \demand^k : k \in \Ntwo }$, and $\bar{\cap}_i = \max\braces{\outages^k_i \cap_i\paren{\time^k} : k \in \Ntwo } $ for all $i \in \gens$. We then define $\bigM_k = \max\braces{D^k - V_k, 0}.$\\

\subsection{Parameter values}
\label{appendix:parameter-values}
\begin{align*}
    \costshed &\in [50000,40000,30000,20000,10000,8000,6500,5000,3000,1500] \tag{Base}\\
    \bound &\in [0, 20, 500, 1000, 1500, 3000, 4500, 7000, 12000, 20000] \tag{BO-$\cvar{}$}\\
    \bound &\in [0, 500, 1000, 2500, 5000, 8000, 10000, 15000, 18000, 30000] \tag{BO-Cond-$\cvar{}$}\\
    \bound &\in [0,10,50,100,250,500,1000,1500,3000,20000] \tag{BO-Cond-$\cvar{}$ (Ind)}
\end{align*}

\subsection{Analysis of impact of temperature dependence in base case study}
\label{app:fig}
\allowbreak
\begin{figure}[H]
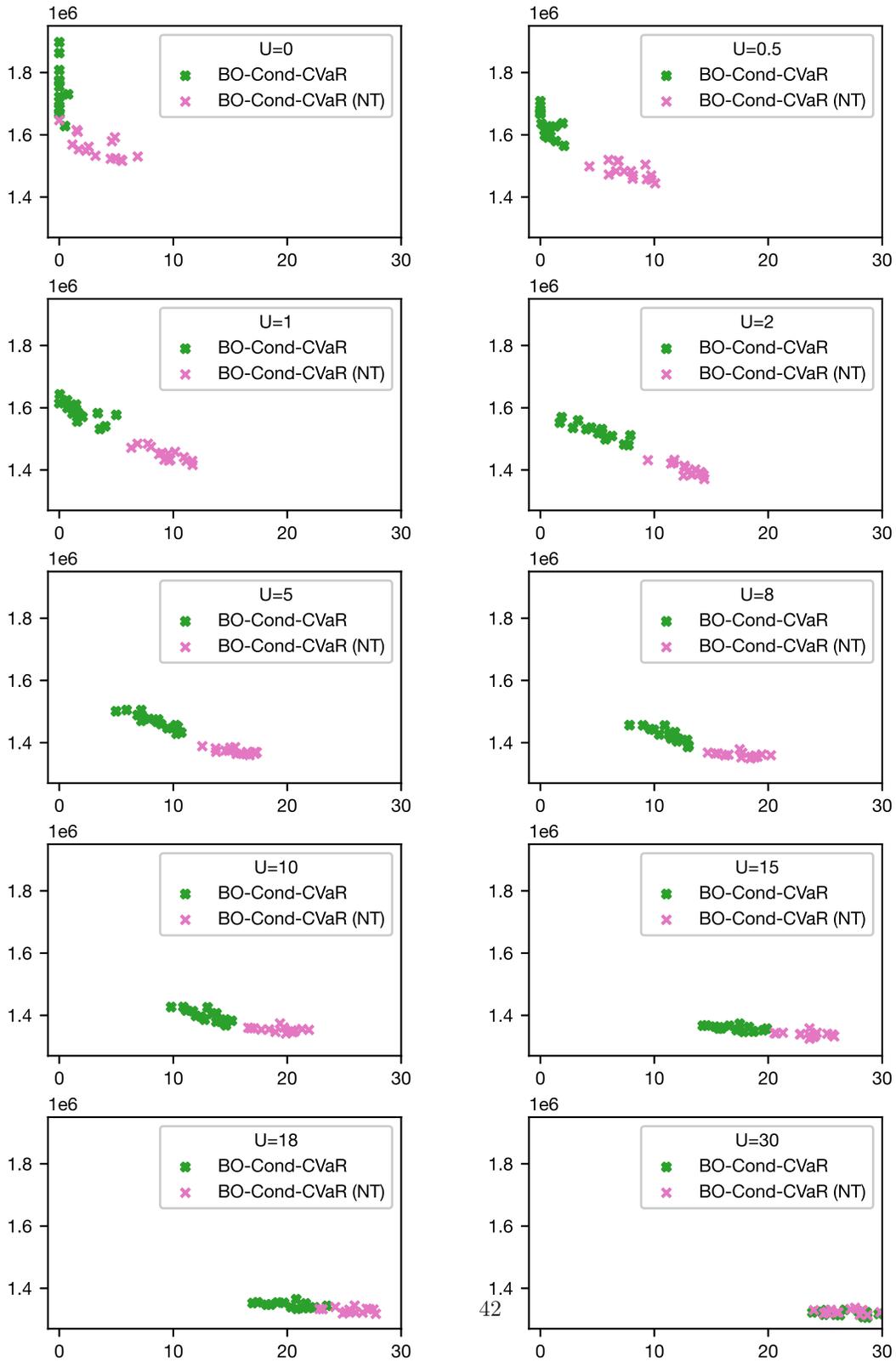


\include{subfigs}

\caption{Effect of generator temperature dependence in base case study by value of second objective upper bound $\bound$. Each figure shows all solutions found for a given value of the $\bound$ for models described in Section \ref{sec:temp-dependence-base}. The bound $\bound$ constrains the conditional load shed ($\cvar{0.1}$, GWh) when solving the models using the $\epsilon$-constraint method. The $x$-axes measure unconditional load shed ($\cvar{0.0001}$, GWh), and the $y$-axes measure total cost. }
\label{fig:notemp-by-U}
\end{figure}

\subsection{Additional generators}
\label{app:gens}
\begin{table}[H]
    \begin{tabular}{ c | c c c c }
Generator & Type & Capacity & Heat rate & Fixed capital cost \\
          &      & (MW)     & (Btu/kWh) & (\$/kW) \\\hline 
CT big & Combustion Turbine & 237 & 9905 & 713 \\
CC small & Combined Cycle & 418 & 6431 & 1084 \\
CC big & Combined Cycle & 1083 & 6370 & 958 \\
Nuclear & Nuclear & 2156 & 10608 & 6041 \\
    \end{tabular}
    \caption{Specifications for new generators added to case study in Section \ref{sec:temp-dependence-modified}.}
    \label{tab:new-generators}
\end{table}

\vspace{0.1cm}
\begin{flushright}
	\scriptsize \framebox{\parbox{2.5in}{Government License: The
			submitted manuscript has been created by UChicago Argonne,
			LLC, Operator of Argonne National Laboratory (``Argonne").
			Argonne, a U.S. Department of Energy Office of Science
			laboratory, is operated under Contract
			No. DE-AC02-06CH11357.  The U.S. Government retains for
			itself, and others acting on its behalf, a paid-up
			nonexclusive, irrevocable worldwide license in said
			article to reproduce, prepare derivative works, distribute
			copies to the public, and perform publicly and display
			publicly, by or on behalf of the Government. The Department of Energy will provide public access to these results of federally sponsored research in accordance with the DOE Public Access Plan. http://energy.gov/downloads/doe-public-access-plan. }}
	\normalsize
\end{flushright}	

\end{document}